\RequirePackage{fix-cm}
\documentclass[smallextended]{svjour3}       
\smartqed  
\usepackage{lineno,hyperref}
\usepackage[fleqn]{amsmath}
\usepackage{amssymb}
\usepackage{latexsym}
\usepackage{epstopdf}
\usepackage{wrapfig}
\usepackage{colortbl}
\usepackage{CJK}
\usepackage{graphicx}
\usepackage{caption}
\usepackage[margin=0.95 in]{geometry}
\begin{document}

\title{Several improved adaptive mapped weighted essentially non-oscillatory scheme for hyperbolic conservation law
}

\author{Shuijiang Tang         
}


\institute{Shuijiang Tang \at
              School of Mathematics and Computational Science, Xiangtan University, Hunan, 411105, China\\              
              \email{sjtang@xtu.edu.cn}           
}

\date{Received: date / Accepted: date}

\maketitle

\begin{abstract}
The decisive factor for the calculation accuracy of the mapped weighted essentially non-oscillatory scheme is the width of the center region of the mapping function. Through analysis of the classical mapped WENO schemes, the results show the width of the central range of the mapping function determined by the local operator in its denominator. Substituting the local operator in WENO-AIM with a symmetric one and an asymmetric function, we get two new adaptive mapped WENO schemes, WENO-AIMS and WENO-AIMA. Similarly, we improve WENO-RM260 and WENO-PM6 by using these local operators, and composed adaptive WENO-RM260 and adaptive WENO-PM6. Theoretical and numerical results show the present adaptive mapped WENO schemes composed in this paper perform better than WENO-AIM, WENO-RM260, and WENO-PM6 for one- and two-dimensional problems.
\end{abstract}

\section{Introduction}
\label{intro}
\par
In recent years, the WENO scheme has been widely used in CFD and has continuously improved. Liu et al. \cite{1} first proposed the idea of the WENO scheme, to assign a reasonable weight to each flux function involved in the essentially non-oscillation (ENO) scheme\cite{2}, and replace it with a weighted convex combination of all fluxes. The basic idea of the WENO scheme is that each flux plays the same role in the smooth region, and only the smoothest flux plays a key role in the region with large gradient changes. In this way, it is possible to maintain the characteristics of substantially essentially no-oscillation and achieve higher accuracy in the smooth region. Jiang and Shu\cite{3} found they can further improve the accuracy of the WENO scheme of Liu et al. And they provided a method to assess the smoothness of templates and a framework for composing the WENO scheme (WENO-JS). Since then, many researchers follow with interest in the WENO scheme and got various related results. In 2000, Balsara and Shu\cite{4} even constructed the WENO scheme with 7 to 13 order precision. In 2003, Qiu et al.\cite{5} got a fifth-order HWENO scheme by using Hermite polynomials. Soon after, they extended the HWENO to two-dimensional problems\cite{6} and the non-uniform grid \cite{7}.
\par
Henrick et al.\cite{8} pointed out the accuracy of the fifth-order WENO-JS scheme will reduce to third-order at critical points in the smooth regions. They reconstructed the weight of the WENO-JS scheme by using a mapping function and got the mapped WENO scheme \raisebox{0.3mm}{------} WENO-M. WENO-M can achieve the optimal convergence order at critical points. After that, Borges et al. \cite{9} proposed an original method to calculate weight by incorporating a global smoothness indicator and constructed the WENO-Z scheme that can achieve the optimal convergence order at critical points. Feng et al. \cite{10} demonstrated the accuracy of the seventh-order or higher-order WENO-M scheme will lower than the WENO-JS scheme in long-term simulation. To address this issue, Feng et al. proposed a new class of piecewise polynomial mapping function by adding two criteria that can decrease the effect from non-smooth stencils and got a new WENO scheme \raisebox{0.3mm}{------} WENO-PM. Compared with WENO-M, WENO-PM6 performs better on both smooth and discontinuous problems. But there are also difficulties in factorization and increased computational cost. Later, Feng et al. \cite{11} gave a new mapping function with two parameters by rewriting the mapping function of the WENO-M method. After selecting suitable parameters, a class of WENO-IM(2,0.1) scheme with less dissipation and high resolution was achieved. Wang et al. \cite{12} gave a new non-segmental mapping function that combined the favorable features of WENO-PM and WENO-IM, and composed a new WENO scheme \raisebox{0.3mm}{------} WENO-RM (260). WENO-RM (260) can increase the non-physical oscillations of seventh and ninth-order WENO schemes. Recently, Vevek et al. \cite{13} adapted the mapping function of WENO-IM with an adaptive operator, developed a rational mapping function with three parameters. They constructed the WENO-AIM ($n, m, s$) scheme, and $n=4,m=2,s=10000$ is subsequently selected. Recently, Hong et al. \cite{14} proposed a pre-discrete mapping method to reduce the cost brought by the mapping process, and given the mapping-function-free WENO-M scheme. Hu \cite{15} applied the mapping function on the indicators of smoothness to increase the efficiency of WENO-M and get a high order WENO-IM scheme.
\par
In this paper, we analyzed the mapping functions of WENO-M, WENO-IM, and WENO-AIM, and found the decisive factor of the three mapping functions are the local operator $\omega_{k}(1-\omega_{k})$ in the denominator. By improving this centrosymmetric local operator, we get a more generalized adaptive mapped WENO scheme \raisebox{0.3mm}{------} WENO-AIMS. We also introduced the asymmetric local operator into the mapping function and getting another adaptive mapped WENO scheme \raisebox{0.3mm}{------} WENO-AIMA. Recommend these local operators in WENO-PM6 and WENO-RM260, we get the adaptive WENO-PM6 and adaptive WENO-RM260 schemes. Many numerical examples show the adaptive WENO-PM and adaptive WENO-RM perform better than WENO-PM and WENO-RM.
 \par
This paper organized as: In Section 2, we briefly introduced the construction of the WENO scheme and several mapping functions. In Section 3, we provide many new mapping WENO schemes and determined suitable parameters. We discuss the performance of the new mapping WENO schemes compared with other mapping WENO schemes for linear advection and Euler equations in Section 4. The conclusions of this paper provided in Section 5.
\section{Numerical method}
In this paper, we study the one-dimensional systems of hyperbolic conservation laws  which is to be solved on the domain $x\in [a,b]$ for $t>0$
\begin{equation}
u_{t}+f(u)_{x}=0,
\end{equation}
where $u(x,t)$ is the conservative variable, $f(u)$ is the flux, $x$ and $t$ are respectively, the space and time.
In the uniform grids, let $x_j=j\triangle x$, $x_{j+ 1/2}=(j+ 1/2)\triangle x,~j=0,1,2,\cdots,N$, where $\triangle x=(b-a)/N$. Eq.(1) at $x_{j}$ can be approximated by
\begin{equation}
\frac{d\bar{u}_j}{dt}=-\frac{1}{\Delta x}(\hat{f}_{j+1/2}-\hat{f}_{j-1/2}),
\end{equation}
where $\hat{f}_{j\pm 1/2}$ are the numerical fluxes.\par
The numerical flux function of the fifth-order WENO scheme can be expressed as:
\begin{equation}
\hat{f}_{j+1/2}=\omega_{0}q_{0}+\omega_{1}q_{1}+\omega_{2}q_{2},
\end{equation}
where $q_{k},~k=0,1,2$ are the third-order flux on the substencil $S_{k}={j-2+k,j-1+k,j+k}$ and given by
\begin{equation}
\begin{cases}
q_{0}=(2f_{j-2}-7f_{j-1}+11f_{j})/6,\\
q_{1}=(-f_{j-1}+5f_{j}+2f_{j+1})/6,\\
q_{2}=(2f_{j}+5f_{j+1}-f_{j+2})/6.
\end{cases}
\end{equation}
\subsection{WENO-JS scheme}
The weights $\omega_{k},~k=0,1,2$ of Jiang and Shu\cite{3} are given as
\begin{equation}
\omega^{JS}_{k}=\frac{\alpha_{k}}{\alpha_{0}+\alpha_{1}+\alpha_{2}}, \alpha_{k}=\frac{d_{k}}{(\beta_{k}+\epsilon)^{2}},~k=0,1,2,
\end{equation}
where $d_{0}=0.1,~d_{1}=0.6$ and $d_{2}=0.3$ are the optimal weights, which generate the fifth-order upstream scheme. $\epsilon$ is a small positive real number that is utilized to avoid the denominator of zero. $\epsilon=10^{-40}$ is recommended by Henrick et al.\cite{8} and used in this paper. And the smoothness indicator $\beta_{k}$ is given by
\begin{equation}
\begin{cases}
\beta_{0}=\frac{13}{12}(f_{j-2}-2f_{j-1}+f_{j})^2+\frac{1}{4}(f_{j-2}-4f_{j-1}+3f_{j})^2,\\
\beta_{1}=\frac{13}{12}(f_{j-1}-2f_{j}+f_{j+1})^2+\frac{1}{4}(f_{j-1}-f_{j+1})^2,\\
\beta_{2}=\frac{13}{12}(f_{j}-2f_{j+1}+f_{j+2})^2+\frac{1}{4}(3f_{j}-4f_{j+1}+f_{j+2})^2.
\end{cases}
\end{equation}
\par
Expanding the smoothness indicators using Taylor series at $x_{j}$ for a smooth solution is given as following
\begin{equation}
\begin{cases}
\beta_{0}=f'^{2}_{j}(\Delta x)^{2}+\left(\frac{13}{12}f_{j}''^{2}-\frac{2}{3}f_{j}'f_{j}'''\right)(\Delta x)^{4}+\left(-\frac{13}{6}f_{j}''f_{j}'''+\frac{1}{2}f_{j}'f_{j}^{(4)}\right)(\Delta x)^{5}+O(\Delta x^{6}),\\
\beta_{1}=f'^{2}_{j}(\Delta x)^{2}+\left(\frac{13}{12}f_{j}''^{2}+\frac{1}{3}f_{j}'f_{j}'''\right)(\Delta x)^{4}+O(\Delta x^{6}),\\
\beta_{2}=f'^{2}_{j}(\Delta x)^{2}+\left(\frac{13}{12}f_{j}''^{2}-\frac{2}{3}f_{j}'f_{j}'''\right)(\Delta x)^{4}+\left(\frac{13}{6}f_{j}''f_{j}'''-\frac{1}{2}f_{j}'f_{j}^{(4)}\right)(\Delta x)^{5}+O(\Delta x^{6}).
\end{cases}
\end{equation}
\par
The expansion (7) means that
\begin{equation}
\beta_{k}=\begin{cases}
f'^{2}_{j}(1+O(\Delta x)^{2}),&f'\neq 0,\\
\frac{13}{12}f''_{j}(\Delta x)^{4}(1+O(\Delta x)),&f'=0.
\end{cases}
\end{equation}
Therefore, WENO-JS suffers from a loss of accuracy near critical points.
\subsection{WENO-M scheme}
The mapped WENO scheme is proposed by Henrick et al\cite{8} by introducing a mapping function. The mapping function is defined as
\begin{equation}
g_{M}(\omega_{k})=\frac{\omega_{k}(d_{k}+d_{k}^{2}-3d_{k}\omega_{k}+\omega_{k}^{2})}{d_{k}^{2}+\omega_{k}(1-2d_{k})},~k=0,1,2.
\end{equation}
\par
The mapping function (9) has following properties:
\begin{equation}
\begin{cases}
monotonically~increasing~in~[0,~1]~with~finite~slopes,\\
g_{M}(0)=0,g_{M}(d_{k})=d_{k},~g_{M}(1)=1,\\
g_{M}'(d_{k})=g_{M}''(d_{k})=0.
\end{cases}
\end{equation}
\par
By using the mapping function, the classical weights mapped to the new ones as follows
\begin{equation}
\omega_{k}^{M}=\frac{g_{M}(\omega_{k})}{g_{M}(\omega_{0})+g_{M}(\omega_{1})+g_{M}(\omega_{2})},~k=0,1,2.
\end{equation}
\par
Applying the Taylor expansion of (9) at $d_{k}$, there is
\begin{equation}
g_{M}(\omega_{k})=d_{k}+O((\Delta x)^{3}).
\end{equation}
\par
Thus, at or near the critical points in smooth regions, the weights of WENO-M satisfies
\begin{equation}
\omega_{k}^{M}=d_{k}+O((\Delta x)^{3}),
\end{equation}
that is, the WENO-M scheme recovers the optimal order of accuracy.
\subsection{WENO-PM scheme}
To improve the accuracy near discontinuities, by adding two criteria, Feng et al. \cite{10} proposed another mapping function
\begin{equation}
g_{PM}(\omega_{k};n)=c_{1}(\omega_{k}-d_{k})^{n+1}(\omega_{k}+c_{2})+d_{k},~n\geq 2,~k=0,1,2,
\end{equation}
where
\begin{equation}
c_{1}=\begin{cases}
(-1)^{n}\frac{n+1}{(d_{k})^{n+1}},~0\leq\omega_{k}\leq d_{k},\\
-\frac{n+1}{(1-d_{k})^{n+1}},~d_{k}\leq\omega_{k}\leq 1,
\end{cases}
c_{2}=\begin{cases}
\frac{d_{k}}{n+1},~0\leq\omega_{k}\leq d_{k},\\
\frac{d_{k}-(n+2)}{n+1},~d_{k}\leq\omega_{k}\leq 1.
\end{cases}
\end{equation}
\par
The mapping function (14) satisfies following properties:
\begin{equation}
\begin{cases}
monotonically~increasing~in~[0,~1]~with~finite~slopes,\\
g_{PM}(0;n)=0,g_{PM}(d_{k};n)=d_{k},~g_{PM}(1;n)=1,\\
g_{PM}'(d_{k};n)=\cdots=g_{PM}^{(n)}(d_{k};n)=0,\\
g_{PM}(1;n)=1,g_{PM}'(1;n)=0.
\end{cases}
\end{equation}
\par
By applying the mapping function and using (16), one can verify that the weights meet
\begin{equation}
\omega^{PM}_{k}=d_{k}+O((\Delta x)^{n+1}),
\end{equation}
this means, WENO-PM scheme can recover the fifth-order accuracy near critical points when $ n\geq 2 $. Feng et al. \cite {10} found that the fifth order WENO-PM method with $n=6$ can get better results than WENO-M. 
\par
Given a tiny positive number $\epsilon$, we call the interval satisfying $\{\omega_{k}|\left|g(\omega_{k})-d_{k}\right|<\epsilon\}$ as a central range of $g(\omega_{k})$. It is clear from Fig. 1 that $g_{PM}(\omega_{k}; 6)$ has a much wider central range compared to $g_{M}(\omega_{k})$. 
\begin{figure}[thb!]
\setlength{\abovecaptionskip}{0pt}
\setlength{\belowcaptionskip}{1pt}
\renewcommand*{\figurename}{Fig.}
\begin{minipage}[t]{0.5\linewidth}
\centering
\includegraphics[scale = 0.5]{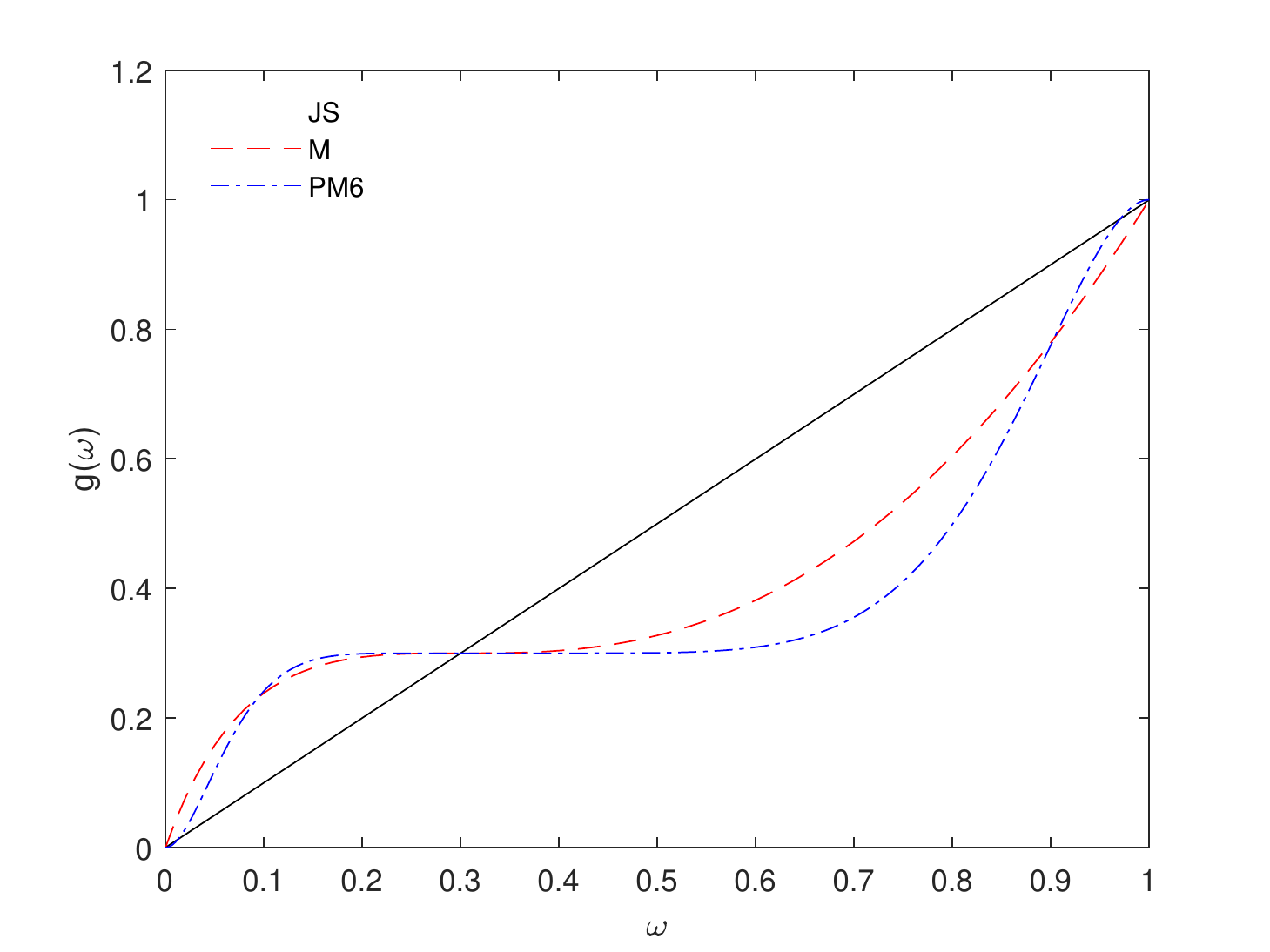}
\label{fig:side:a}
\small
\caption{Comparison of mapping functions for $d=0.3$}
\end{minipage}%
\begin{minipage}[t]{0.5\linewidth}
\centering
 \includegraphics[scale = 0.5]{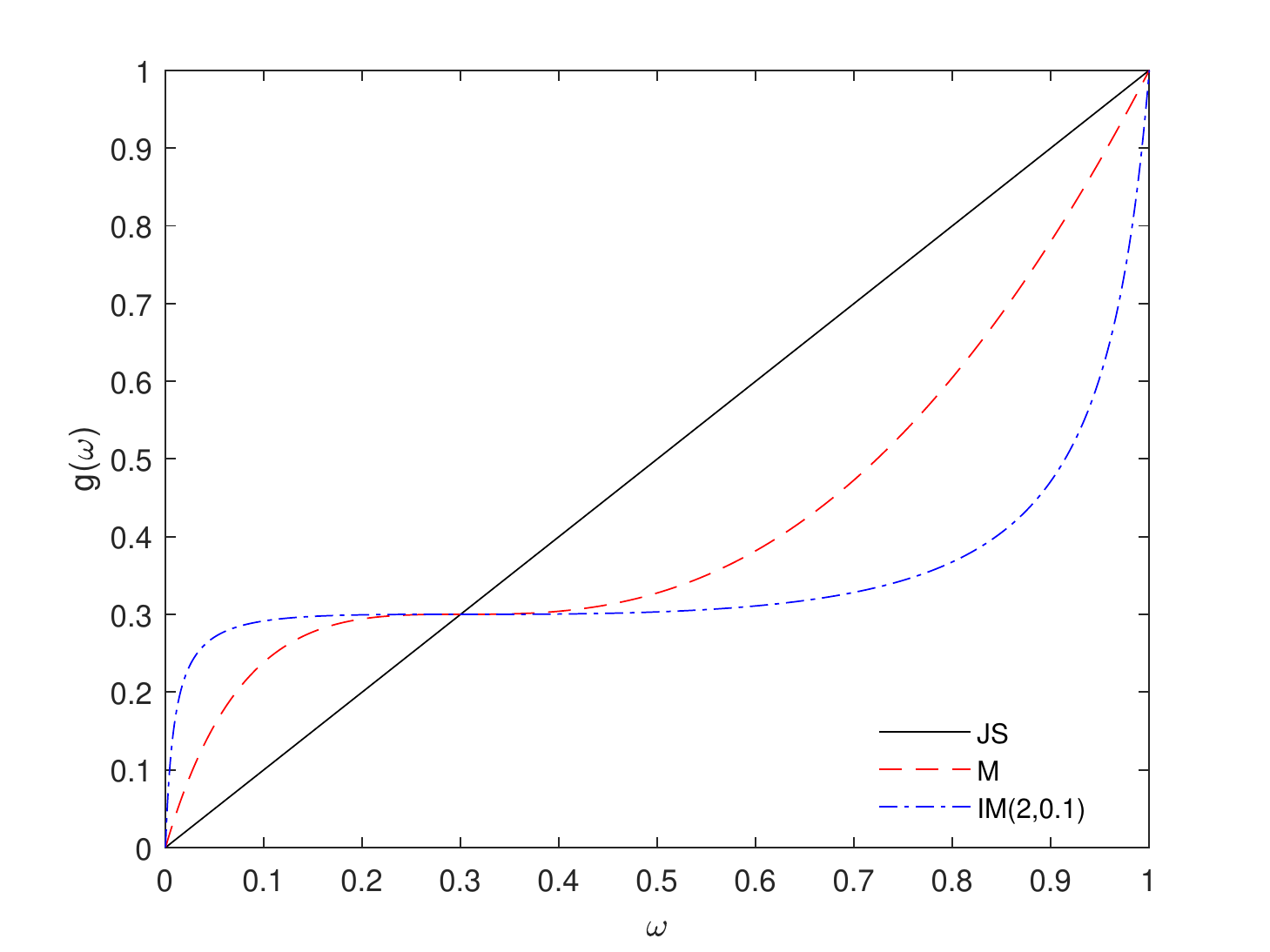}
\label{fig:side:b}
\small
\caption{Comparison of mapping functions for $d=0.3$}
\end{minipage}
\end{figure}
\subsection{WENO-IM scheme}
Feng et al. \cite{11} developed an improved version of the WENO-M method; i.e., WENO-IM. By rewriting and changing the original mapping function, they got a new type of mapping functions with two parameters, $A$ and $n$, as follows:
\begin{equation}
g_{IM}(\omega_{k};n,A)=d_{k}+\frac{(\omega_{k}-d_{k})^{n+1}A}{(\omega_{k}-d_{k})^{m}A+\omega_{k}(1-\omega_{k})},
\end{equation}
where $A>0$, $n$ is a positive even integer. The mapping function (18) also satisfies the property of (16).
Applying the Taylor expansion of (18) at $d_{k}$, there is
\begin{equation}
g_{IM}(\omega_{k};n,A)=d_{k}+O((\Delta x)^{n+1}).
\end{equation}
\par
Feng et al.\cite{11} found that by taking $n=2$ and $A = 0.1$, WENO-IM performs better than WENO-M. It is evident from Fig. 2 that $g_{IM}(\omega_{k}; 2,0.1)$ has a much wider central range compared to $g_{M}(\omega_{k})$.
\subsection{WENO-RM scheme}
Wang et al\cite{12} developed an improved version of the WENO-M method, i.e., WENO-RM.  Reconstructed the mapping function of WENO-M, they obtained a new mapping function with three parameters, $m$, $n$ and $\tau$, as follows
\begin{equation}
g_{RM}(\omega_{k};m,n,\tau)=d_{k}+\frac{(\omega_{k}-d_{k})^{n+1}}{a_{0}+a_{1}\omega_{k}+\cdots+a_{m+s+1}\omega_{k}^{m+\tau+1}},
\end{equation}
where $m\leq n$, and the choice of coefficients $a_{i}$ is related to the value of $\tau$, usually there are
\begin{equation}
1)~\tau=0,\begin{cases}
a_{i}=C_{n+1}^{i}(-d_{k})^{n-i},~i=0,1,2,\cdots,m,\\
a_{m+1}=(1-d_{k})^{n}-\displaystyle\sum_{i=0}^{m}a_{i}.
\end{cases}
\end{equation}
\begin{equation}
2)~\tau=1,\begin{cases}
a_{i}=C_{n+1}^{i}(-d_{k})^{n-i},~i=0,1,2,\cdots,m,\\
a_{m+1}=(m+2)\left[(1-d_{k})^{n}-\displaystyle \sum_{i=0}^{m}a_{i}\right]+\sum_{i=0}^{m}ia_{i}-(n+1)(1-d_{k})^{n-1},\\
a_{m+2}=(n+1)(1-d_{k})^{n-1}-\displaystyle \sum_{i=0}^{m}ia_{i}-(m+1)\left[(1-d_{k})^{n}-\sum_{i=0}^{m}a_{i}\right].
\end{cases}
\end{equation}
\par
The mapping function (20) satisfies following properties:
\begin{equation}
\begin{cases}
monotonically~increasing~in~[0,~1]~with~finite~slopes,\\
g_{RM}(0;m,n,\tau)=0,g'_{RM}(0;m,n,\tau)=\cdots=g^{(m)}_{RM}(0;m,n,\tau)=0,\\
g_{RM}(d_{k};m,n,\tau)=d_{k},g'_{RM}(d_{k};m,n,\tau)=\cdots=g^{(n)}_{RM}(d_{k};m,n,\tau)=0,\\
g_{RM}(1;m,n,\tau)=1,g'_{RM}(1;m,n,\tau)=\cdots=g^{(\tau)}_{RM}(1;m,n,\tau)=0.
\end{cases}
\end{equation}
\par
Applying the Taylor expansion of (20) at $d_{k}$, there is
\begin{equation}
g_{RM}(\omega_{k};m,n,\tau)=d_{k}+O((\Delta x)^{n+1}).
\end{equation}
\par
Wang et al.\cite{12} found that,  taking $m=2$, $n=6$ and $\tau=0$, the  fifth and seventh order  WENO-RM perform better than  WENO-M. It is evident from Fig. 3 that $g_{RM}(\omega_{k}; 2,6,0)$ has a much wider central range compared to $g_{M}(\omega_{k})$.
\begin{figure}[thb!]
\setlength{\abovecaptionskip}{0pt}
\setlength{\belowcaptionskip}{1pt}
\renewcommand*{\figurename}{Fig.}
\begin{minipage}[t]{0.5\linewidth}
\centering
\includegraphics[scale = 0.5]{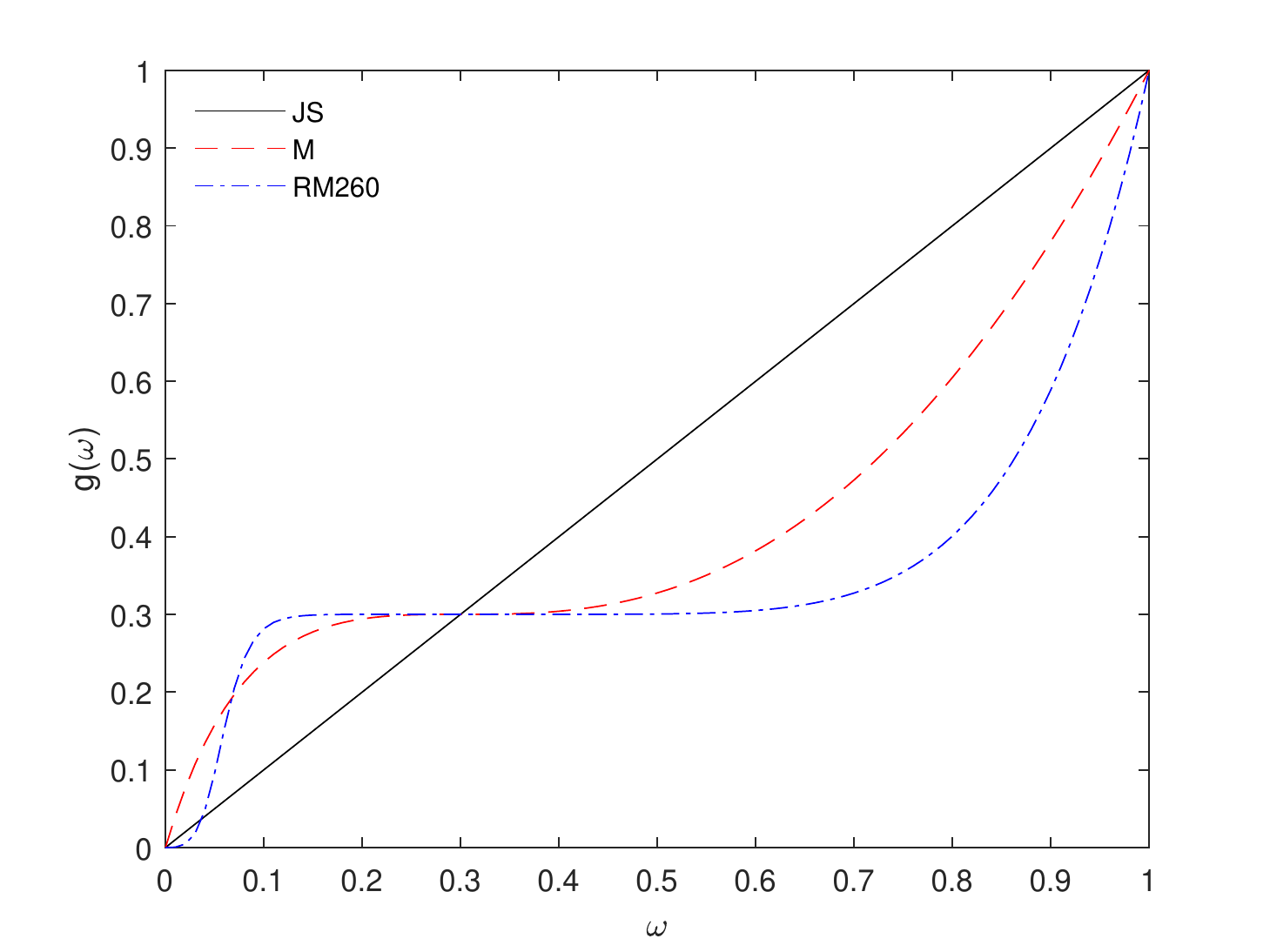}
\label{fig:side:a}
\small
\caption{Comparison of mapping functions for $d=0.3$}
\end{minipage}%
\begin{minipage}[t]{0.5\linewidth}
\centering
 \includegraphics[scale = 0.5]{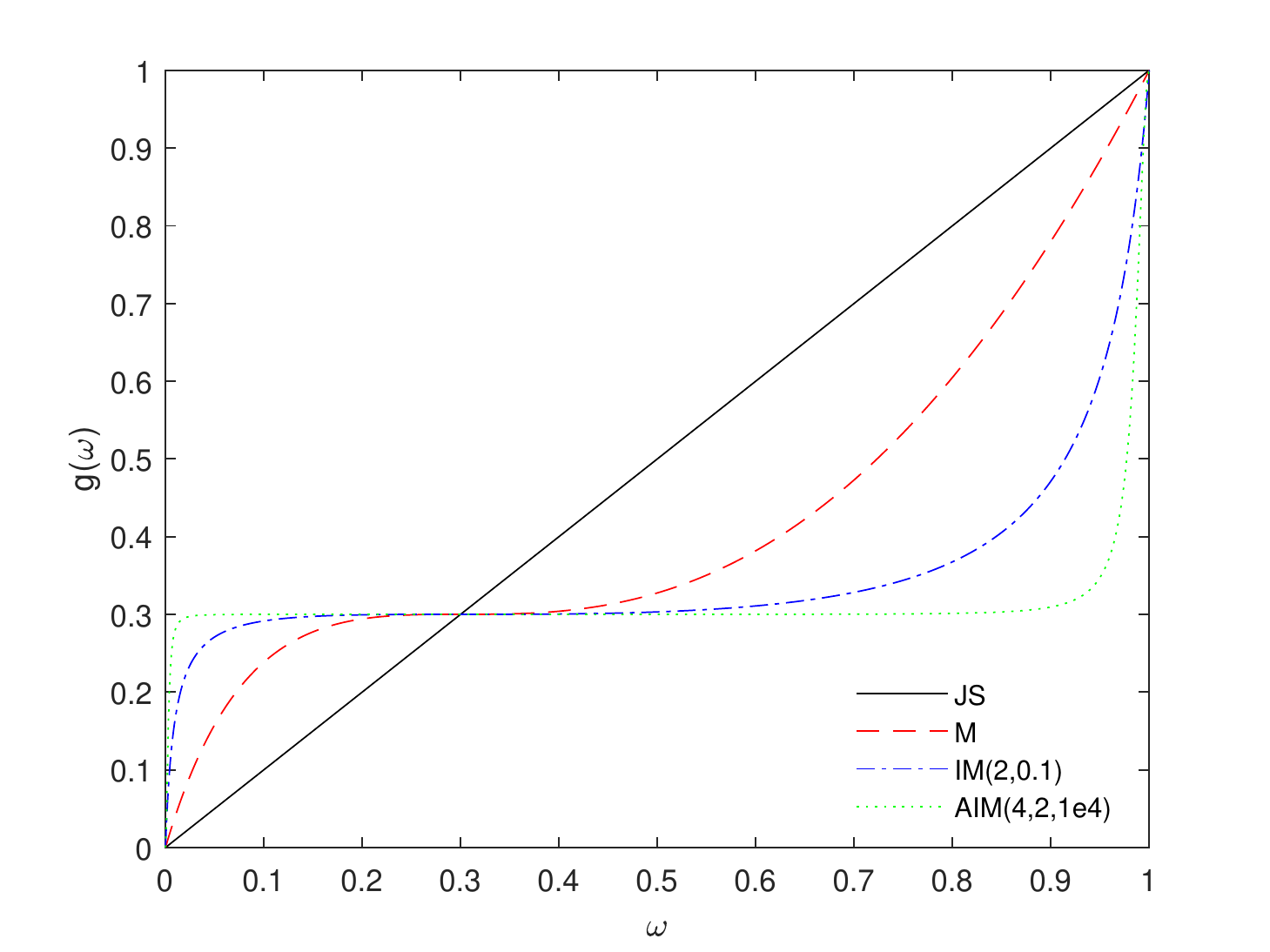}
\label{fig:side:b}
\small
\caption{Comparison of mapping functions for $d=0.3$}
\end{minipage}
\end{figure}
\subsection{WENO-AIM scheme}
Recently, Vevek et al.\cite{13} improved the mapping function of WENO-IM, developed a rational mapping function with three parameters, and constructed the adaptive WENO-AIM ($n, m, s$) scheme.
\begin{equation}
g_{AIM}(\omega_{k};n,m,s)=d_{k}+\frac{(\omega_{k}-d_{k})^{n+1}}{(\omega_{k}-d_{k})^{n}+s(\omega_{k}(1-\omega_{k}))^{m}},
\end{equation}
where $n$ is positive even number, $m$ positive  number and $s=cd_{k}^{-1}\lambda$ with
\begin{equation}
\lambda=\frac{min(\beta_{k})}{max(\beta_{k})+\epsilon_{m}},
\end{equation}
where, $\epsilon_{m}$ is a small number that is provides a threshold of significance for the smoothness indicators, and usually selected $\epsilon_{m}=(\Delta x)^{5}$.
\par
The mapping function (25) satisfies following properties:
\begin{equation}
\begin{cases}
monotonically~increasing~in~[0,~1]~when~n\geq m-1,\\
g_{AIM}(0;n,m,s)=0,~g_{AIM}(d_{k};n,m,s)=d_{k},~g_{AIM}(1;n,m,s)=1,\\
g'_{AIM}(0;n,m,s)=\cdots=g^{(n)}_{AIM}(0;n,m,s)=0,\\
g'_{AIM}(0;n,m,s)=g'_{AIM}(1;n,m,s)=1~for~m>1.
\end{cases}
\end{equation}
\par
Applying the Taylor expansion of (25) at $d_{k}$, there is
\begin{equation}
g_{AIM}(\omega_{k};n,m,s)=d_{k}+O((\Delta x)^{n+1}).
\end{equation}
\par
Vevek et al.\cite{13} found that by taking $m=2$, $n=4$ and $c=10^{4}$, far better results could be obtained for the seventh order WENO methods  compared to WENO-M and WENO-IM(2,0.1). It is evident from Fig. 4 that $g_{AIM}(\omega_{k}; 4,2,s=1e4)$ has a much wider central range compared to $g_{M}(\omega_{k})$ and $g_{IM}(\omega_{k};2,0.1)$.
\section{Design and properties of new mapped WENO schemes}
\subsection{The improvement of WENO-AIM scheme}
From the mapping functions of WENO-M, WENO-IM and WENO-AIM, there exists a common formula in these three formulas
\begin{equation}
\phi(\omega_{k})=(\omega_{k}(1-\omega_{k}))^{\kappa},
\end{equation}
it is the determining factor in determining the center width of the mapping function, and we call it the width operator in this paper. And it meets the following conditions
\begin{equation}
\begin{cases}
\phi~is~a~continuous~positive~function~on~[0,~1], \\
\phi~is~monotone~increase~in~[0,~\frac{1}{2}],~and~decrease~in~[\frac{1}{2},~1],\\
\phi(0)=\phi(1)=0,\\
\phi(\omega_{k})=\phi(1-\omega_{k}),\\
\phi'(0)=\phi'(1)=0,~for~\kappa>1.
\end{cases}
\end{equation}
\par
As a matter of fact, we can always improve the above conditions (30) as
\begin{equation}
\begin{cases}
\phi~is~a~continuous~positive~function~on~[0,~1], \\
\phi(0)=\phi(1)=\phi'(0)=\phi'(1)=0.\\
There ~exist~\omega_{L},~\omega_{R}\in (0,1),~\omega_{L}\leq\omega_{R},~\phi~increases~in~(0,~\omega_{L})~and~decreases~in~(\omega_{R},~1).
\end{cases}
\end{equation}
\par
Therefore, we can extend the mapping function of WENO-AIM to the following form
\begin{equation}
g_{AIM\phi}(\omega_{k};n,s,\phi)=d_{k}+\frac{(\omega_{k}-d_{k})^{n+1}}{(\omega_{k}-d_{k})^{n}+s\phi(\omega_{k})},
\end{equation}
\par
Differentiate the above formula with respect to $\omega_{k}$, we have
\begin{equation}
g'_{AIM\phi}(\omega_{k};n,s,\phi)=\frac{(\omega_{k}-d_{k})^{n}\left[(\omega_{k}-d_{k})^{n}+s((n+1)\phi(\omega_{k})-\phi'(\omega_{k})(\omega_{k}-d_{k}))\right]}{\left[(\omega_{k}-d_{k})^{n}+s\phi(\omega_{k})\right]^{2}}.
\end{equation}
\par
Obviously, when the local operator $\phi(\omega_{k})$ satisfies all conditions in (31), and
\begin{equation}
(n+1)\phi(\omega_{k})-\phi'(\omega_{k})(\omega_{k}-d_{k})\geq 0,
\end{equation}
then, $g_{AIM\phi}(\omega_{k};n,s,\phi)$ satisfies the properties of (27). In other words, we can always build any local operator that supports the implementation (31) and (34) so that the function constructed by (32) becomes a mapping function that satisfies the condition (27). The resulting adaptive mapped WENO scheme is called the generalized WENO-AIM scheme.
\par
Next, we will use the different locations of the $\omega_{L},~\omega_{R}$ in (31) to create two new types of local operators that are different from those of AIM.
\subsubsection{Center symmetric local operator}
When $\omega_{L}=1-\omega_{R}$, and the width function $\phi(\omega_{k})$ meets the conditions (31), (33) and
\begin{equation}
\phi(\omega_{k})=\phi(1-\omega_{k}),
\end{equation}
it is called the center symmetric local operator. Obviously, the local operator in AIM is of this kind. In addition, the symmetric local operator can also be developed in the following form.
\begin{equation}
~\phi(\omega_{k})=(1+\chi(\omega_{k}-0.5)^{2})(\omega_{k}(1-\omega_{k}))^{\kappa}.
\end{equation}
where $\kappa$ is an integer greater or equal to 2, $\chi$ is undetermined non-negative constant. Easy to verify, the new local operator satisfy the conditions of (31) and (34), and when $\chi=0$, it is the local operator of the WENO-AIM scheme. Figure 5 shows the graph of the function $\phi(\omega_{k})$ as $\kappa=2$. One can see from the graph, the values of $\phi(\omega_{k})$ increase with the increase of $\chi$. Figure 6 shows the mapping weight function constructed by using (36) with the effect of varying parameter $\chi$, where $n=4$, $\kappa=2$ and $s=1$. From this figure, the width of the central range becomes wider as $\chi$ increases. This shows, the larger $\chi$ is, the closer the weight composed by the mapping function $g_{AIM\phi}(\omega_{k};n,s,\phi)$ is to the optimal weight.
\begin{figure}[thb!]
\setlength{\abovecaptionskip}{0pt}
\setlength{\belowcaptionskip}{1pt}
\renewcommand*{\figurename}{Fig.}
\begin{minipage}[t]{0.5\linewidth}
\centering
\includegraphics[scale = 0.52]{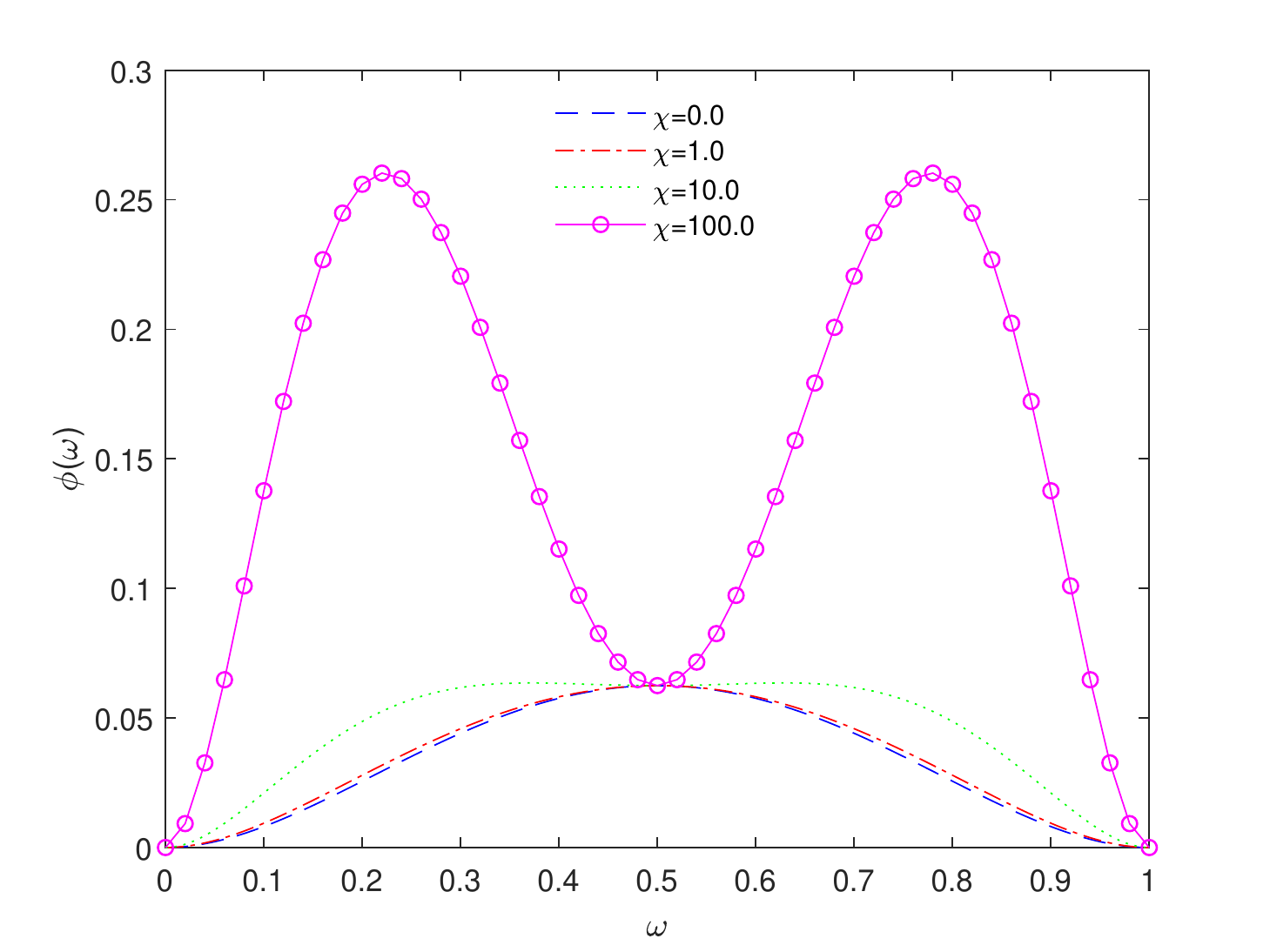}
\label{fig:side:a}
\small
\caption{The curve of $\phi(\omega_{k})$ for (36) with $\kappa=2$}
\end{minipage}%
\begin{minipage}[t]{0.5\linewidth}
\centering
 \includegraphics[scale = 0.52]{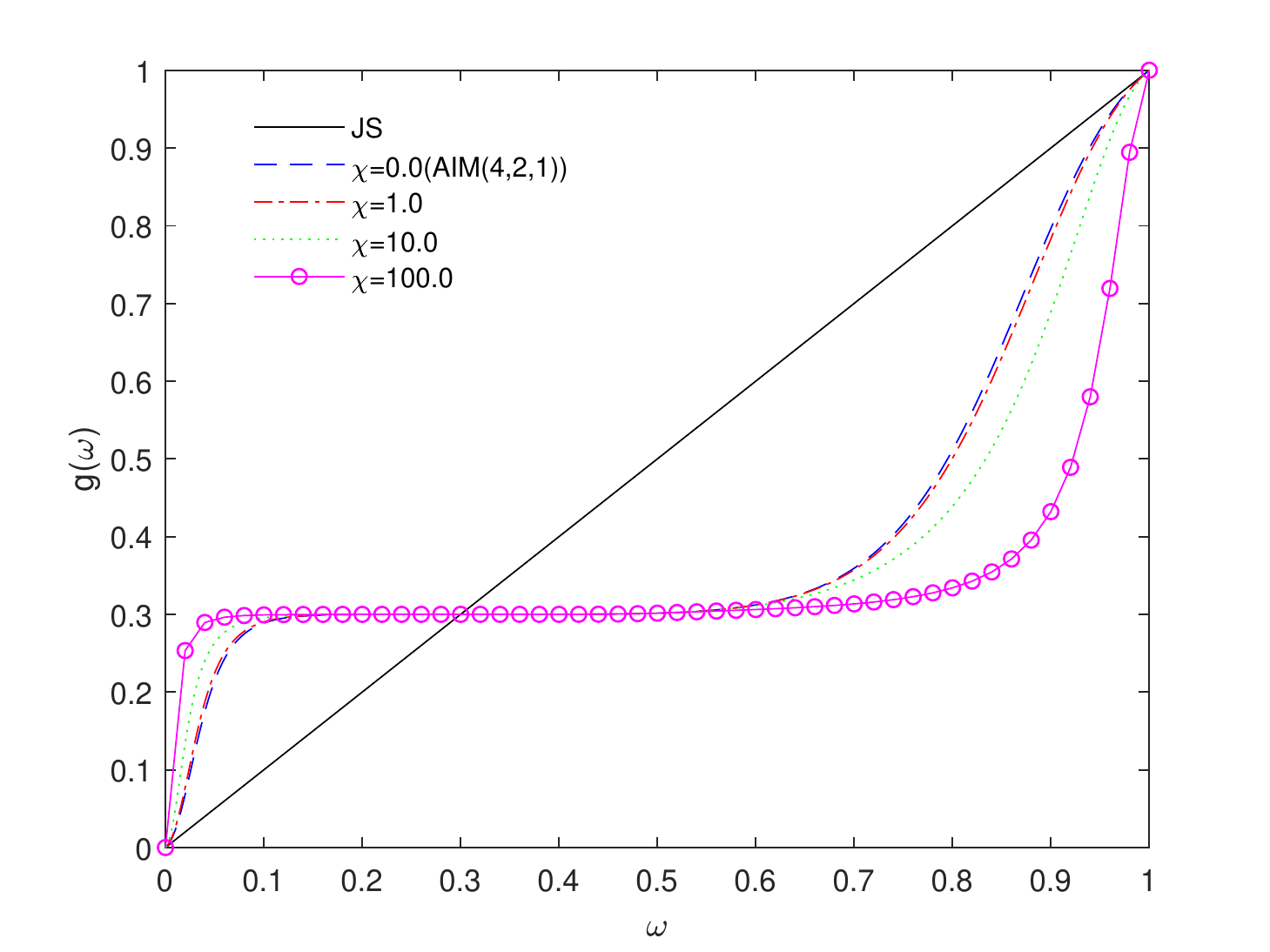}
\label{fig:side:b}
\small
\caption{Effect of varying parameter $\chi$ on $g_{AIM\phi}(\omega_{k};n,s,\phi)$
for (36) at $n=4$, $\kappa=2$ and $s=1$ }
\end{minipage}
\end{figure}
\par
\textbf{Remark 1}~In addition to (36), we can also construct several center symmetric local operators in the following form
\begin{equation}
\begin{aligned}
~\phi(\omega_{k})&=P_{m}(\omega_{k})(\omega_{k}(1-\omega_{k}))^{\kappa},\\
~\phi(\omega_{k})&=min(\omega_{k}^{\kappa},(1-\omega_{k})^{\kappa}).\\
~\phi(\omega_{k})&=sin(\pi\omega_{k})^{\kappa},...
\end{aligned}
\end{equation}
where $\kappa$ is is an integer greater or equal to 2, $P_{m}(\omega_{k})$ is a positive polynomial of degree $m$, and satisfies $P_{m}(\omega_{k})=P_{m}(1-\omega_{k})$.
\subsubsection{Asymmetric local operator}
Corresponding to the symmetric local operator constructed above, we can also construct an asymmetric local operator in the following form
\begin{equation}
~\phi(\omega_{k})=(1+\chi\omega_{k})(\omega_{k}(1-\omega_{k}))^{\kappa}.
\end{equation}
where $\kappa$ is an integer greater than or equal to 2, and $\chi$ is an undetermined constant, with $\chi>-1$. Easy to verify, the new local operator meets the requirements of (31) and (34), and then is the WENO-AIM scheme as $\chi=0$. Figure 7 shows the graph of the function $\phi(\omega_{k})$ with the variable parameter $\chi$. As seen from the graph, the values of $\phi(\omega_{k})$ increase by $\chi$, and then the critical point gradually moves to the right. We also compared $g_{AIM}(\omega_{k};4,2,1)$ and $g_{AIM\phi}(\omega_{k};4,1,\phi)$ to (38) for $\kappa=2$, as shown in Figure 8. One can see that the width of the central range of $g_{AIM\phi}(\omega_{k};4,1,\phi)$ is wider than $g_{AIM}(\omega_{k};4,2,1)$ as $\chi>0$. But we also noticed that the central part gradually moves to the right as $\chi$ increases.
\par
\textbf{Remark 2}~In addition to (38), we can also construct several asymmetric local operators as following
\begin{equation}
\begin{aligned}
~\phi(\omega_{k})&=P_{m}(\omega_{k})(\omega_{k}(1-\omega_{k}))^{\kappa},\\
~\phi(\omega_{k})&=min(\omega_{k}^{\kappa},(1-\omega_{k})^{\iota}),...
\end{aligned}
\end{equation}
where $\kappa$ and $\iota$ are unequal positive integers greater than 2, $P_{m}(\omega_{k})$ is a positive polynomial of degree $m$, and satisfies $P_{m}(\omega_{k})\neq P_{m}(1-\omega_{k})$.
\par
Usually, we call the $AIM\phi$ uses the center symmetric local operator (36) or (37) as AIMS, and the mapped WENO scheme uses AIMS \raisebox{0.3mm}{------} WENO-AIMS. The $AIM\phi$ uses the asymmetric local operator (38) or (39) is called AIMA, and the mapped WENO scheme uses AIMA is called WENO-AIMA. In this paper, we will only consider the two scenarios (36) and (38).
\begin{figure}[htp]
\setlength{\abovecaptionskip}{0pt}
\setlength{\belowcaptionskip}{1pt}
\renewcommand*{\figurename}{Fig.}
\begin{minipage}[t]{0.5\linewidth}
\centering
\includegraphics[scale = 0.5]{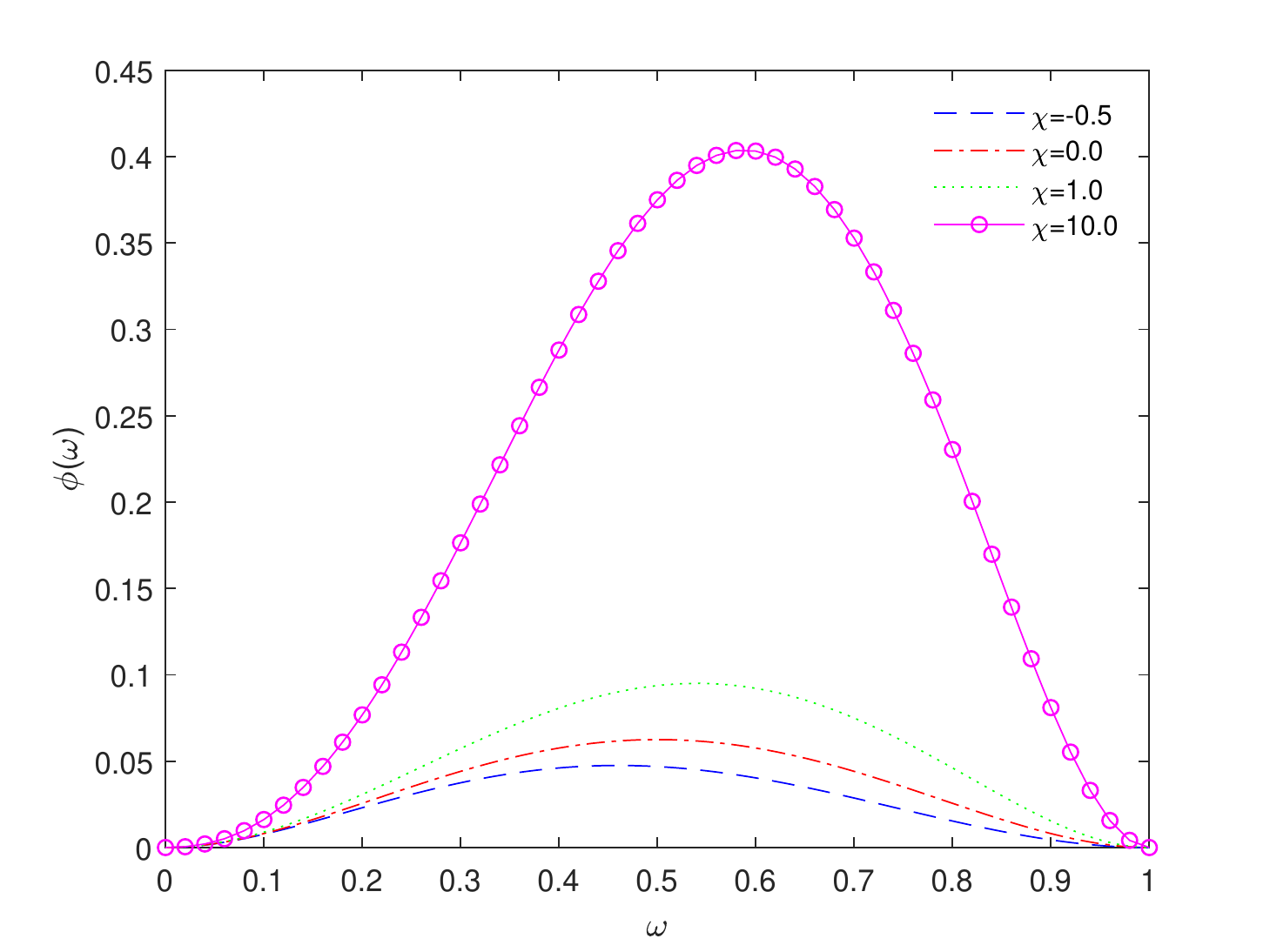}
\label{fig:side:a}
\small
\caption{The curve of $\phi(\omega_{k})$ for (37) with $\kappa=2$}
\end{minipage}%
\begin{minipage}[t]{0.5\linewidth}
\centering
 \includegraphics[scale = 0.5]{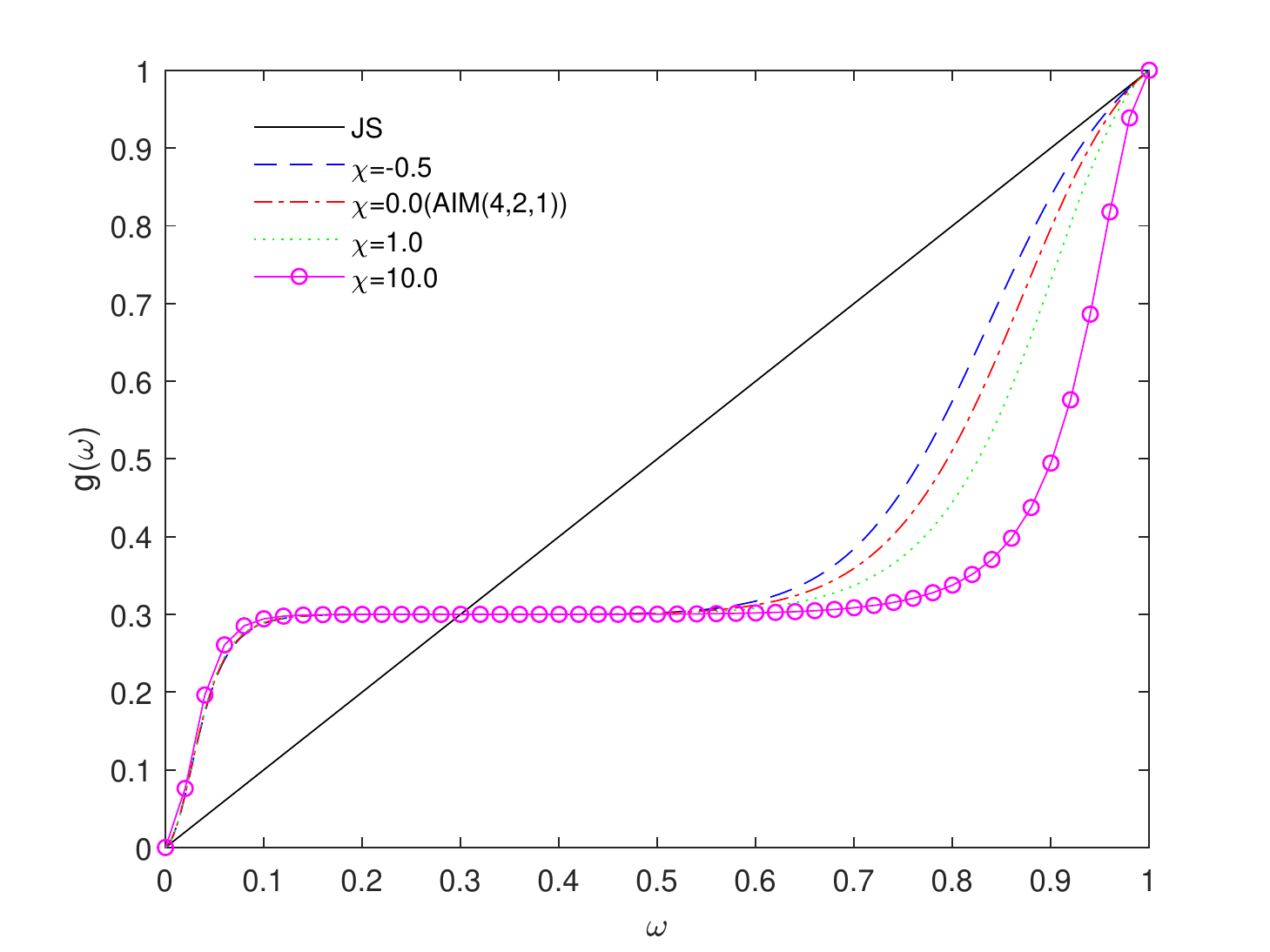}
\label{fig:side:b}
\small
\caption{Effect of varying parameter $\chi$ on $g_{AIM\phi}(\omega_{k};n,s,\phi)$
for (37) at $n=4$, $\kappa=2$ and $s=1$}
\end{minipage}
\end{figure}
\subsection{The adaptive WENO-PM scheme}
Although WENO-PM can achieve the optimal accuracy at the critical points, it will be unstable during long-term calculations. Therefore, we will improve it by using the previously local operator. We can always rewrite the mapping function of WENO-PM into the following form
\begin{equation}
g_{PM}(\omega_{k};n)=d_{k}+\frac{(\omega_{k}-d_{k})^{n+1}}{\frac{1}{c_{1}(\omega_{k}+c_{2})}}.
\end{equation}
\par
Using the construction idea of $AIM\phi$, we can get an adaptive PM mapping function of the form
\begin{equation}
g_{APM}(\omega_{k};n,s,\phi)=d_{k}+\frac{(\omega_{k}-d_{k})^{n+1}}{\frac{1}{c_{1}(\omega_{k}+c_{2})}+s\phi(\omega_{k})},
\end{equation}
where $n$ is a positive even number, $c_{1}$ and $c_{2}$ defined as (19), $\phi(\omega_{k})$ defined as (36) and (37), and $s$ defined as
\begin{equation}
s=\frac{cd_{k}min(\beta_{k})}{max(\beta_{k})+(\triangle x)^{5}}.
\end{equation}
\par
Easy to verify, the adaptive PM mapping function satisfies the conditions of (18). When the solution is smooth, $s=cd_{k}$, so as $c$ tends to infinity, there is $g_{APM}(\omega_{k};n,s,\phi)=d_{k}$. When the solution is discontinuous, there is $g_{APM}(\omega_{k};n,s,\phi)=g_{PM}(\omega_{k};n)$. The improved WENO-PM scheme using (41) will always achieve the optimal accuracy when the solution is continuous. And it will degenerate to WENO-PM when the solution is discontinuous. We call the improved WENO-PM scheme equips with the symmetric local operator (36) as WENO-APMS, and the improved WENO-PM scheme equips with the asymmetric local operator (38) as WENO-APMA.
\par
Figure 9 shows the curve of the mapping function $g_{PM}$, $g_{APMS}$ and $g_{APMA}$ with $n=6$, $s=10$ and $\chi=1$. One can see from the graph, the width of the central range of $g_{APMS}$ and $g_{APMA}$ is the same and wider than $g_{PM}$. This shows that when the values of $s$ and $\chi$ are the same, the width of the center region of the mapping function $g_{APMS}$ and $g_{APMA}$ is the same, but we can also see that $g_{APMA}$ is closer to the optimal weight than $g_{APMS}$.
\begin{figure}[thb!]
\setlength{\abovecaptionskip}{0pt}
\setlength{\belowcaptionskip}{1pt}
\renewcommand*{\figurename}{Fig.}
\begin{minipage}[t]{0.5\linewidth}
\centering
\includegraphics[scale = 0.5]{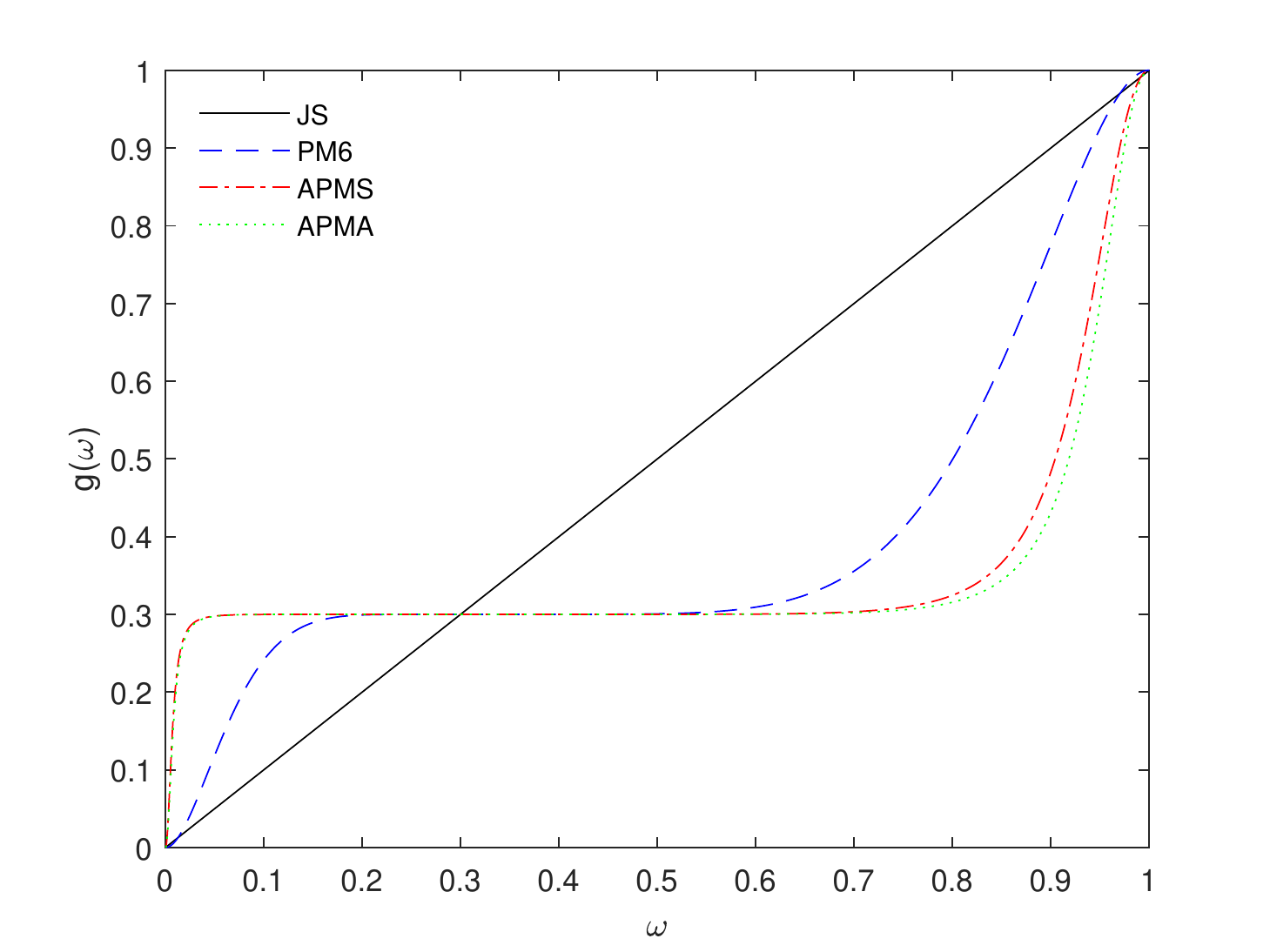}
\label{fig:side:a}
\small
\caption{Comparison of $g_{PM}$, $g_{APMS}$  and $g_{APMA}$ with $n=2$, \protect\\
$s=10$ and $\chi=1$.}
\end{minipage}%
\begin{minipage}[t]{0.5\linewidth}
\centering
 \includegraphics[scale = 0.5]{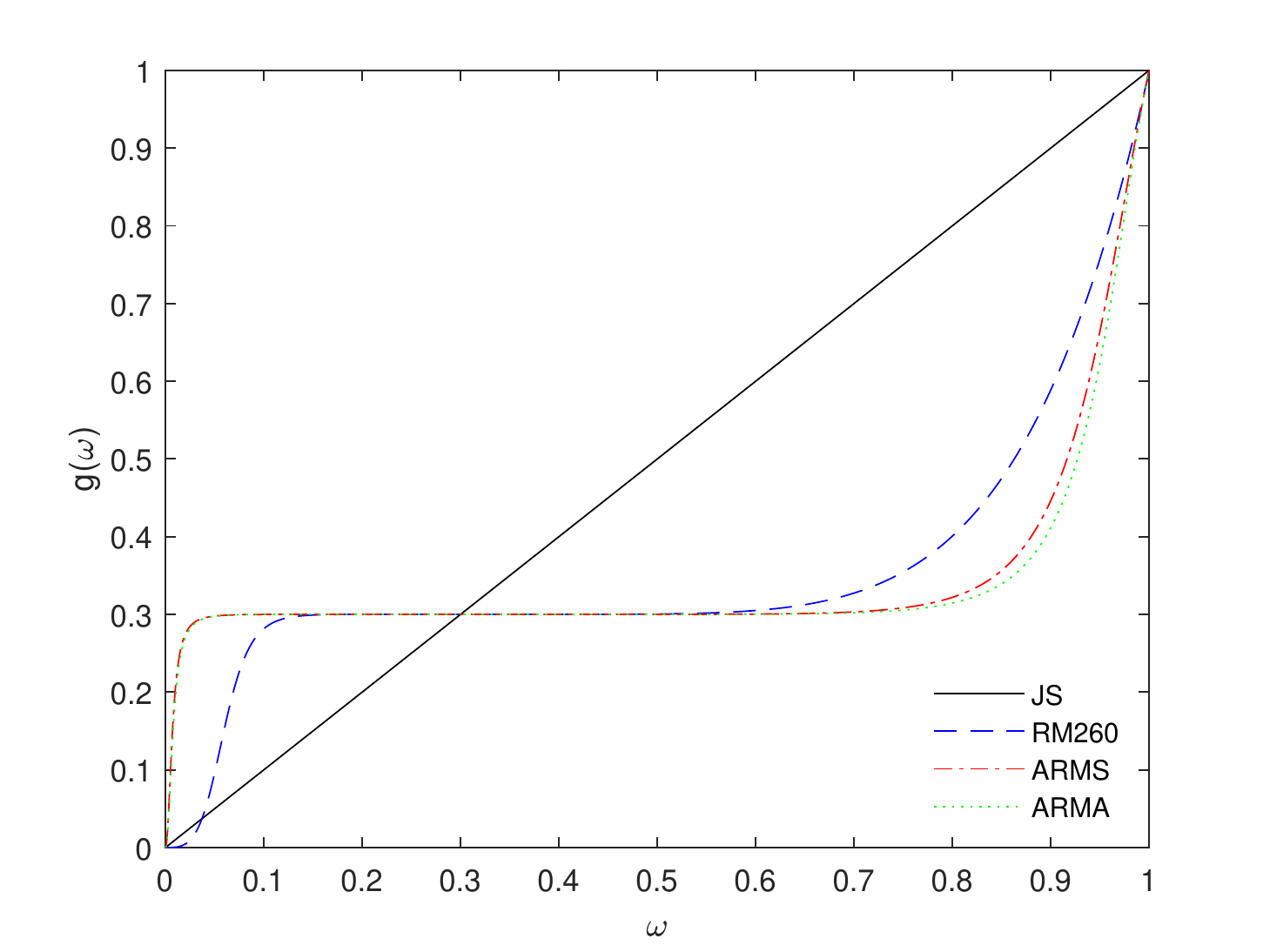}
\label{fig:side:b}
\small
\caption{Comparison of $g_{RM}$, $g_{ARMS}$ and $g_{ARMA}$ with $m=2,$ \protect\\
$n=6,~\tau=0$, $s=10$ and $\chi=1$}
\end{minipage}
\end{figure}
\par
\textbf{Remark 3}~ Noted that the definition of $s$ in (42) differs from the definition of AIM. If use $s$ in AIM, the new scheme constructed in this section will oscillate in the calculation of various problems, which will lead to instability.
\subsection{The adaptive WENO-RM scheme}
Using the construction idea of $AIM\phi$, we can get an adaptive RM mapping function of the form
\begin{equation}
g_{ARM}(\omega_{k};m,n,\tau,s,\phi)=d_{k}+\frac{(\omega_{k}-d_{k})^{n+1}}{a_{0}+a_{1}\omega_{k}+\cdots+a_{m+\tau+1}\omega_{k}^{m+\tau+1}+s\phi(\omega_{k})},
\end{equation}
where $m$ and $n$ are positive numbers, $m\leq n$, and the choice of coefficients $a_{i}$ defined as (21),(22), $s$ defined as (42), $\phi(\omega_{k})$ defined as (36) and (38). Easy to verify, the adaptive RM mapping function satisfies the conditions of (23). We call the improved WENO-RM scheme equips with (36) as WENO-ARMS, and the improved WENO-RM scheme equips with (38) as WENO-ARMA.
\par
Figure 10 shows the curve of the mapping function $g_{RM}$ and the adaptive mapping function $g_{ARMS}$ and $g_{ARMA}$ with $m=2,~n=6,~\tau=0$, $s=10$ and $\chi=1$. The width of the central range of both $g_{ARMS}$ and $g_{ARMA}$ is the same and wider than $g_{RM}$. This shows that when the values of $s$ and $\chi$ are the same, the width of the center region of the mapping function $g_{ARMS}$ and $g_{ARMA}$ is the same, and the value of $g_{APMA}$ is closer to the optimal weight than $g_{APMS}$.
\subsection{Spectral properties of new mapped WENO schemes}
Using the approximate dispersion relation of shock-capturing schemes (ADR) provided by \cite{20}, we first analyze the spectral properties of WENO-AIMS and WENO-AIMA. Although there are two parameters to determine in WENO-AIMS and WENO-AIMA. But whether it is WENO-AIMS or WENO-AIMA, it is the extension of WENO-AIM. Thus, we always take the value of $c$ in WENO-AIMS and WENO-AIMA as 10000, and the spectral resolution of both schemes for selected $\chi$ are given in Figs. 11 and 12.
\par
One can see from Fig. 11, the dispersion and dissipation curves of WENO-AIMS will get closer to the fifth-order upwind scheme with the increase of $\chi$. In Figure 12, the variation trend of the dispersion curve is consistent with the value of $\chi$. The dissipation curve when $\chi=1$ is the closest to the fifth upwind scheme and is almost the same as $\chi=10,100$ and $1000$. In general, the WENO-AIMS and WENO-AIMA with positive $\chi$ have better spectral resolution in all wave numbers compared with that of WENO-AIM.
\begin{figure}[thb!]
\setlength{\abovecaptionskip}{0pt}
\setlength{\belowcaptionskip}{0pt}
\renewcommand*{\figurename}{Fig.}
\centering
\includegraphics[scale = 0.66]{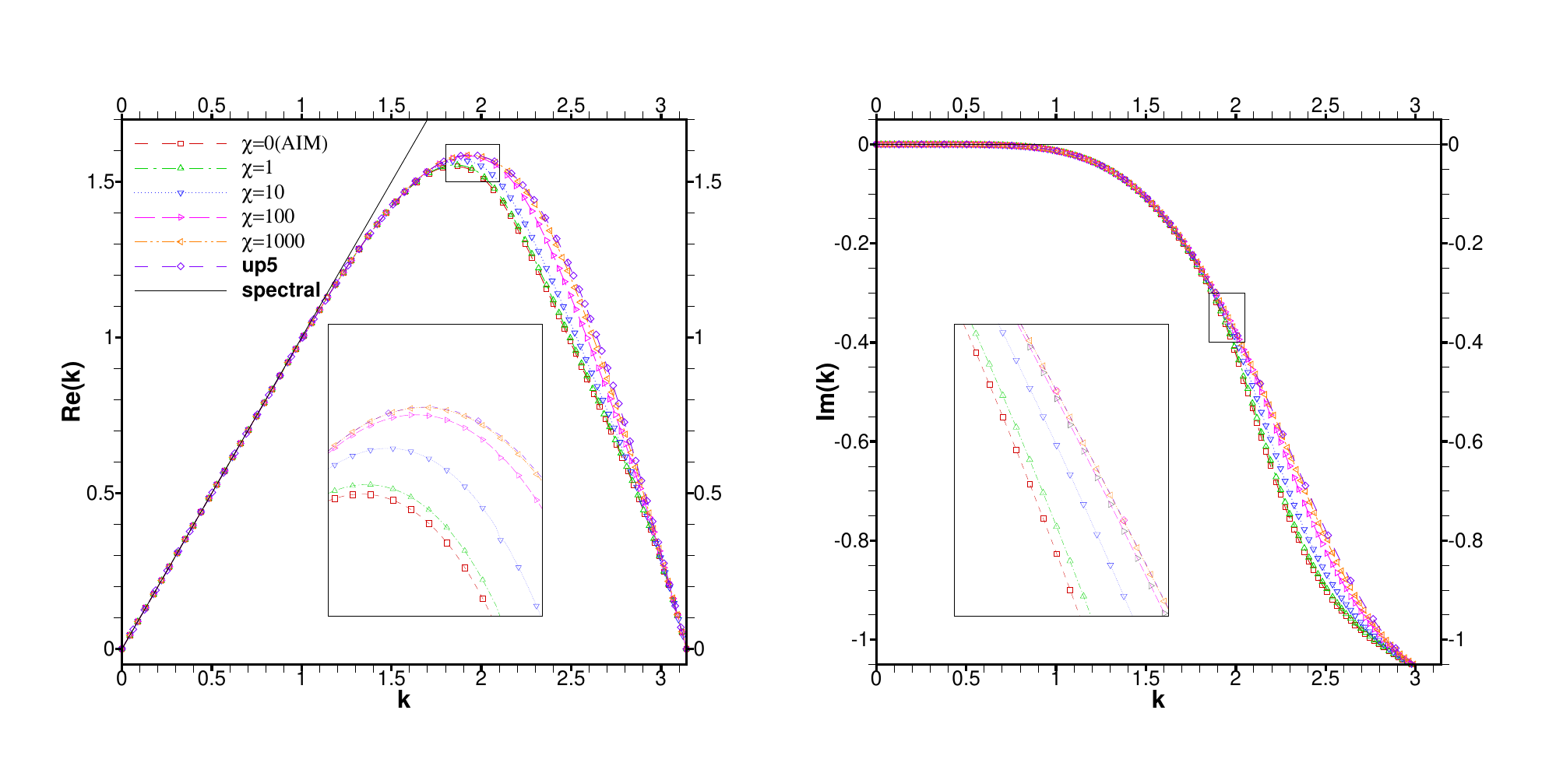}
\small
\caption{Dispersion (left) and dissipation (right) of the WENO-AIMS scheme with selected $\chi$ by the ADR analysis.}
\end{figure}
\begin{figure}[thb!]
\setlength{\abovecaptionskip}{0pt}
\setlength{\belowcaptionskip}{0pt}
\renewcommand*{\figurename}{Fig.}
\centering
\includegraphics[scale = 0.66]{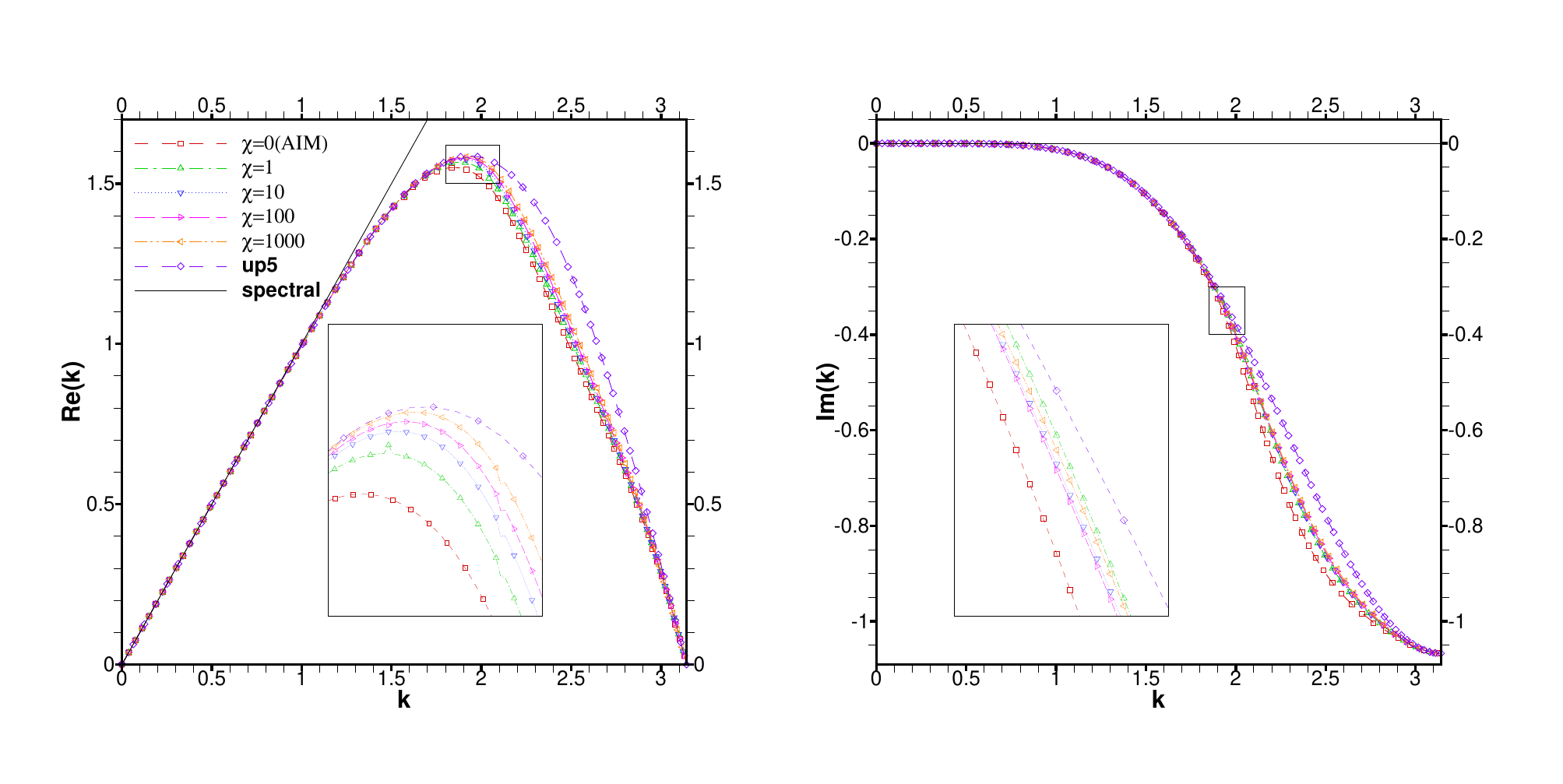}
\small
\caption{Dispersion (left) and dissipation (right) of the WENO-AIMA scheme with selected $\chi$ by the ADR analysis.}
\end{figure}
\begin{figure}[thb!]
\setlength{\abovecaptionskip}{0pt}
\setlength{\belowcaptionskip}{0pt}
\renewcommand*{\figurename}{Fig.}
\centering
\includegraphics[scale = 0.66]{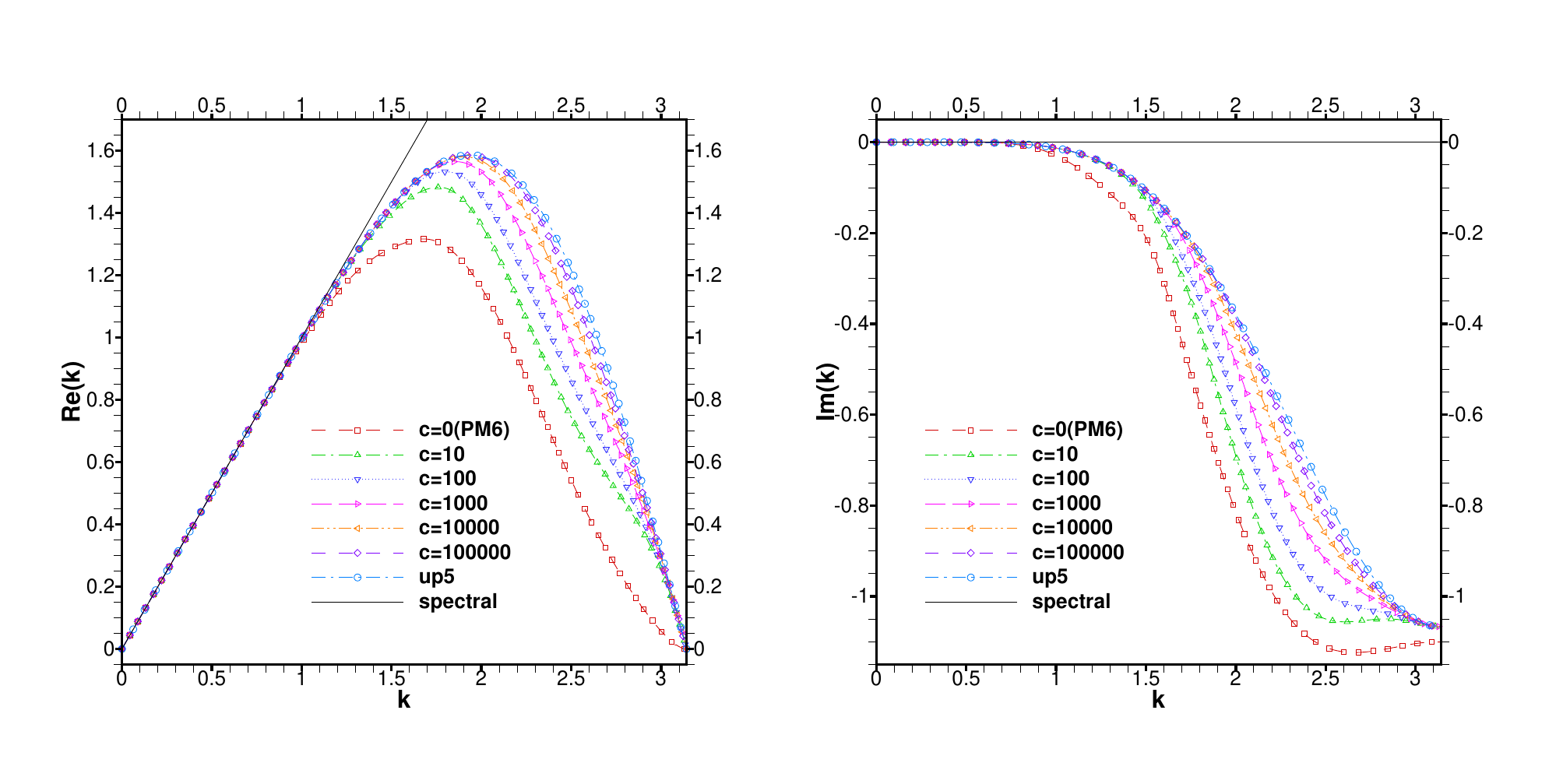}
\small
\caption{Dispersion (left) and dissipation (right) of the WENO-APMS scheme with selected $c$ when $\chi=100$ by the ADR analysis.}
\end{figure}
\begin{figure}[thb!]
\setlength{\abovecaptionskip}{0pt}
\setlength{\belowcaptionskip}{0pt}
\renewcommand*{\figurename}{Fig.}
\centering
\includegraphics[scale = 0.66]{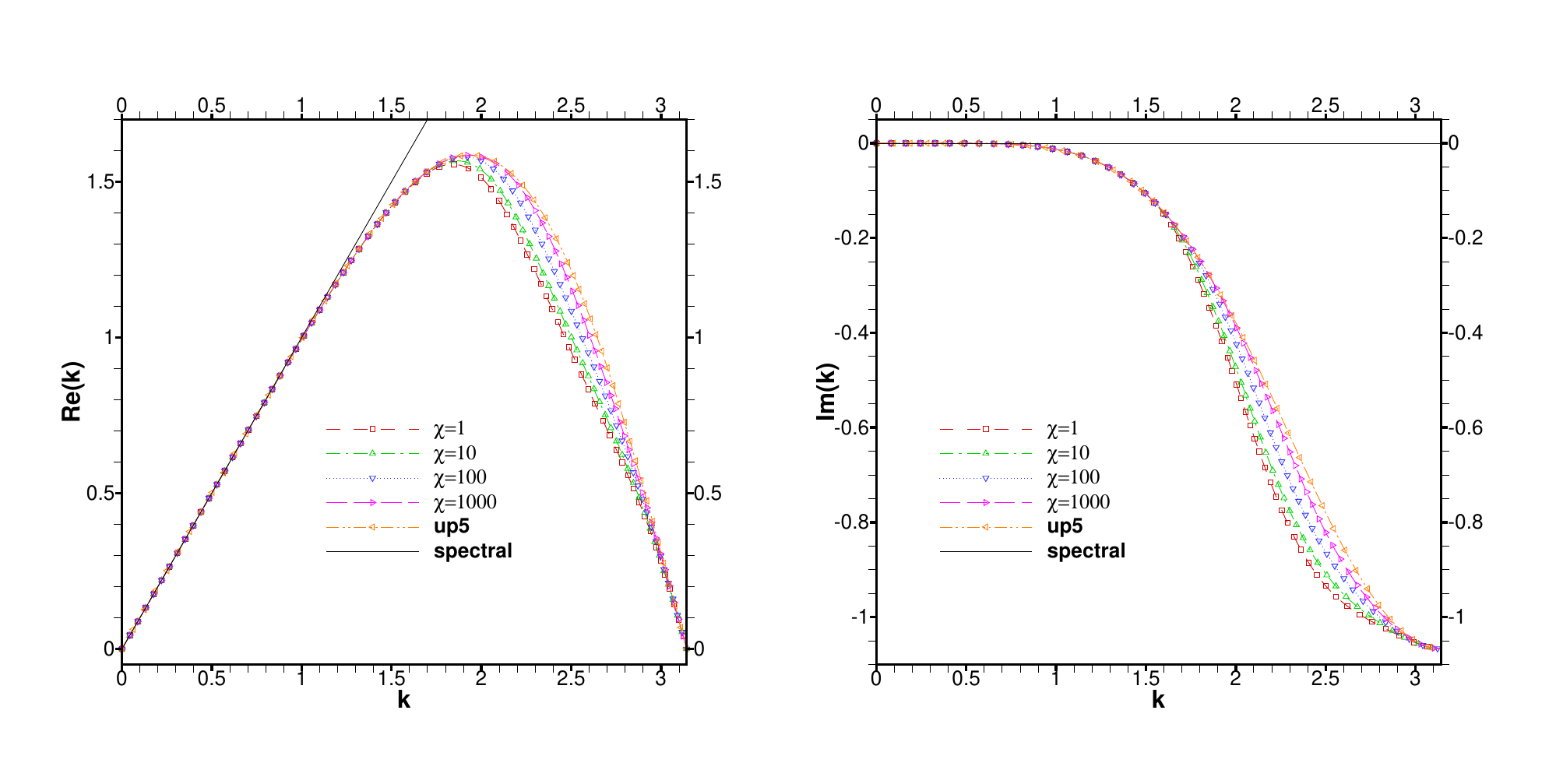}
\small
\caption{Dispersion (left) and dissipation (right) of the WENO-APMS scheme with selected $\chi$ when $c=10000$ by the ADR analysis.}
\end{figure}
\begin{figure}[thb!]
\setlength{\abovecaptionskip}{0pt}
\setlength{\belowcaptionskip}{0pt}
\renewcommand*{\figurename}{Fig.}
\centering
\includegraphics[scale = 0.66]{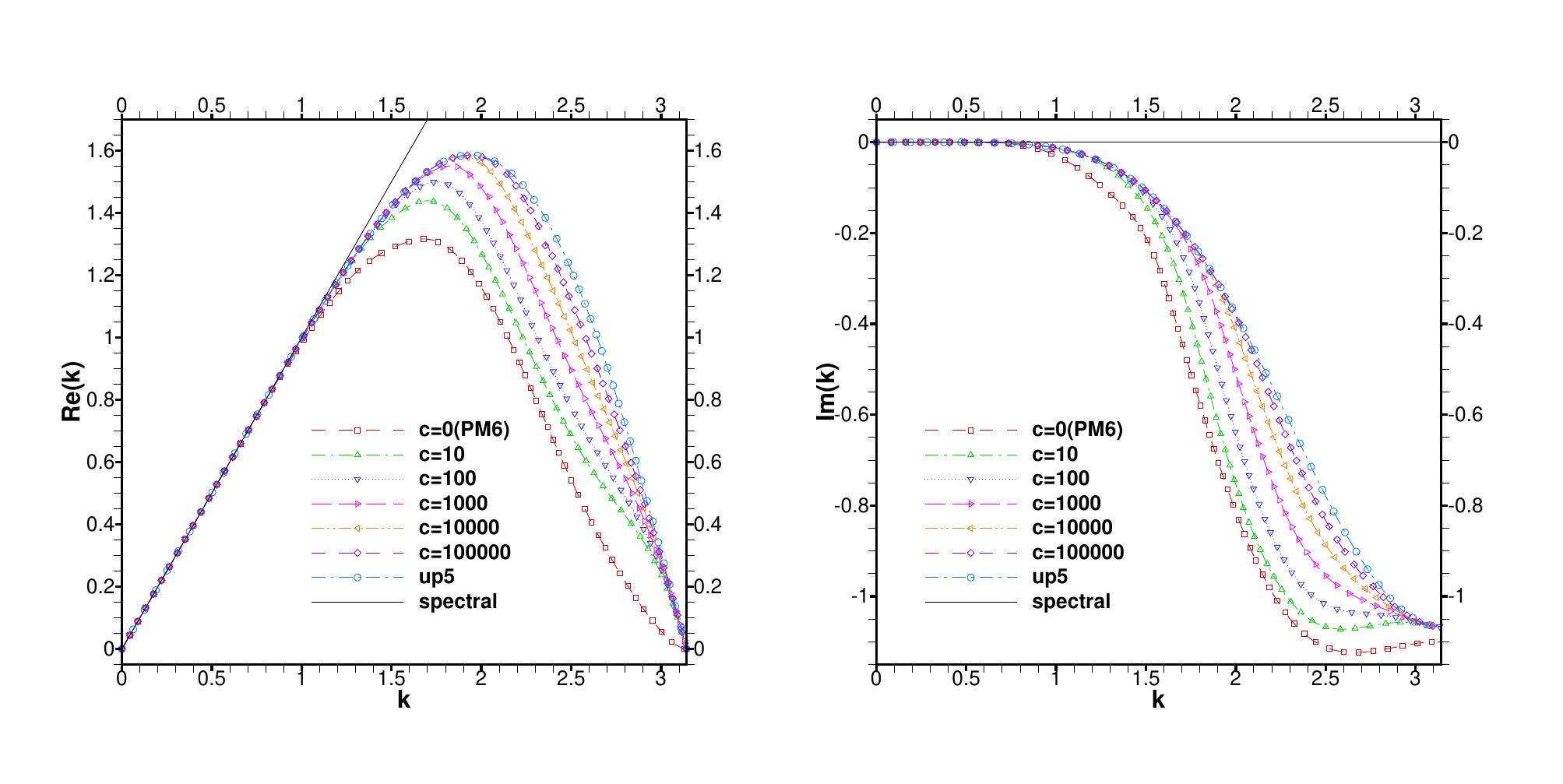}
\small
\caption{Dispersion (left) and dissipation (right) of the WENO-APMA scheme with selected $c$ when $\chi=100$ by the ADR analysis.}
\end{figure}
\begin{figure}[thb!]
\setlength{\abovecaptionskip}{0pt}
\setlength{\belowcaptionskip}{0pt}
\renewcommand*{\figurename}{Fig.}
\centering
\includegraphics[scale = 0.66]{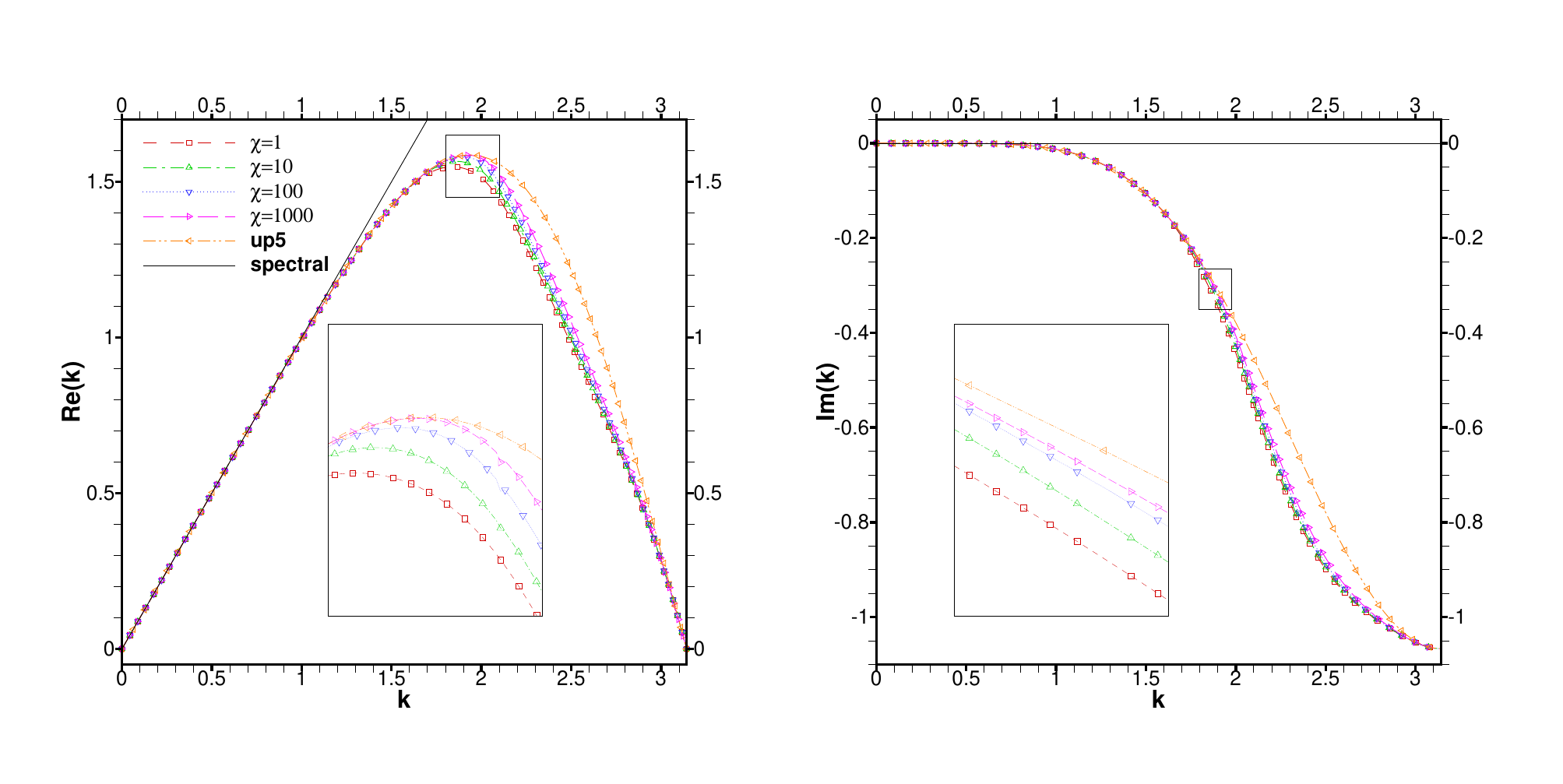}
\small
\caption{Dispersion (left) and dissipation (right) of the WENO-APMA scheme with selected $\chi$ when $c=10000$ by the ADR analysis.}
\end{figure}
\begin{figure}[thb!]
\setlength{\abovecaptionskip}{0pt}
\setlength{\belowcaptionskip}{0pt}
\renewcommand*{\figurename}{Fig.}
\centering
\includegraphics[scale = 0.66]{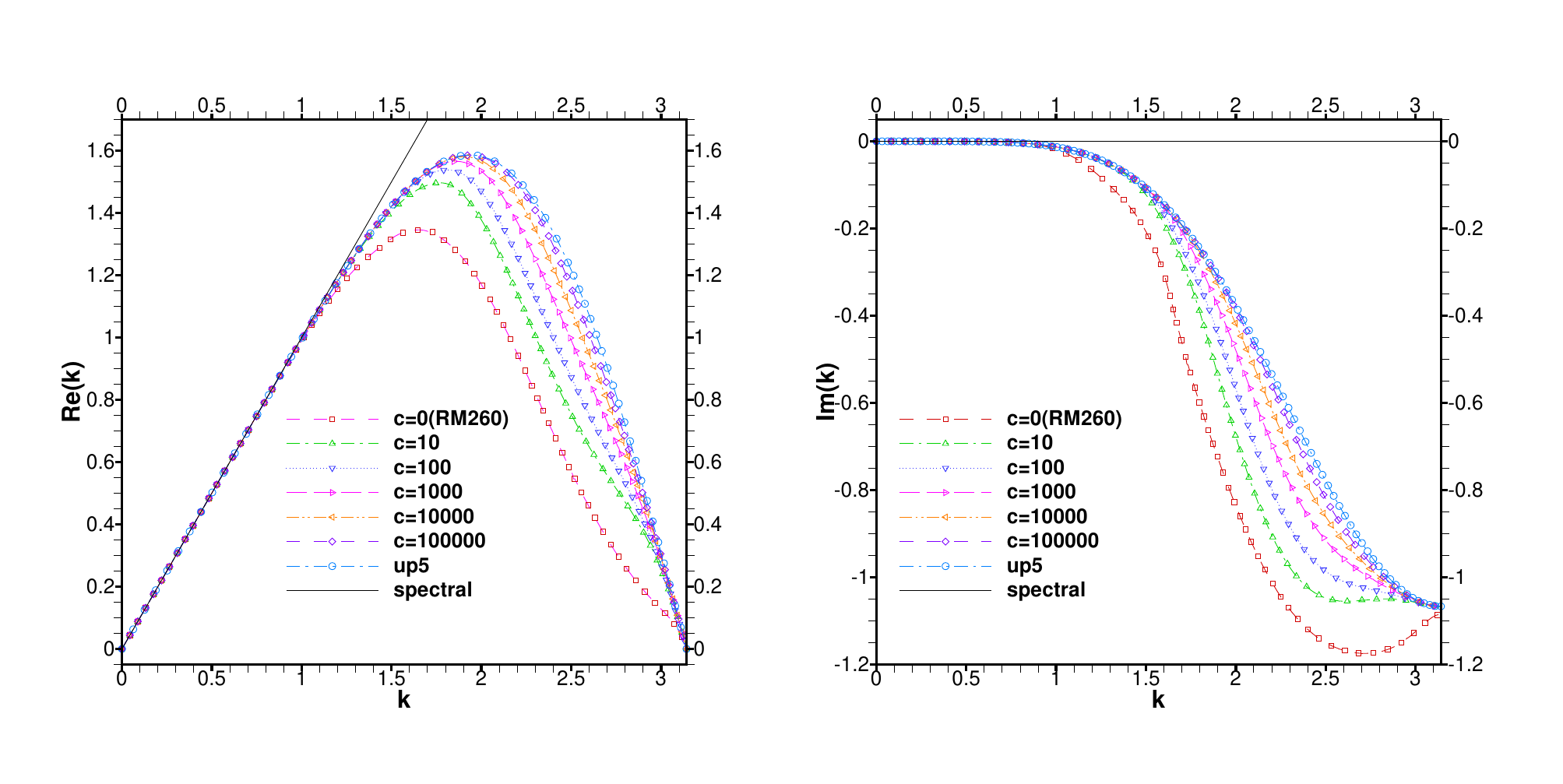}
\small
\caption{Dispersion (left) and dissipation (right) of the WENO-ARMS scheme with selected $c$ when $\chi=100$ by the ADR analysis.}
\end{figure}
\begin{figure}[thb!]
\setlength{\abovecaptionskip}{0pt}
\setlength{\belowcaptionskip}{0pt}
\renewcommand*{\figurename}{Fig.}
\centering
\includegraphics[scale = 0.66]{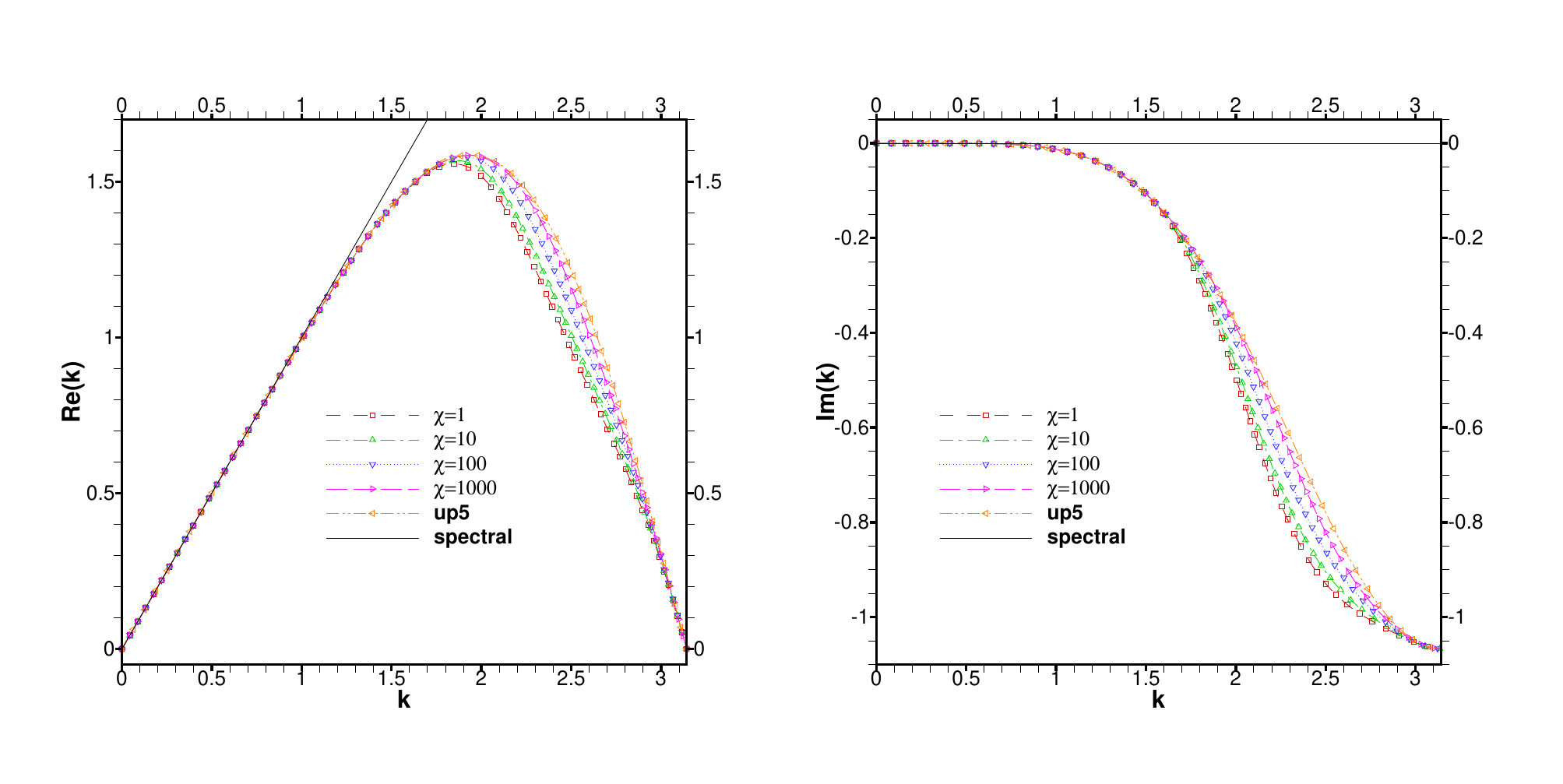}
\small
\caption{Dispersion (left) and dissipation (right) of the WENO-ARMS scheme with selected $\chi$ when $c=10000$ by the ADR analysis.}
\end{figure}
\begin{figure}[thb!]
\setlength{\abovecaptionskip}{0pt}
\setlength{\belowcaptionskip}{0pt}
\renewcommand*{\figurename}{Fig.}
\centering
\includegraphics[scale = 0.66]{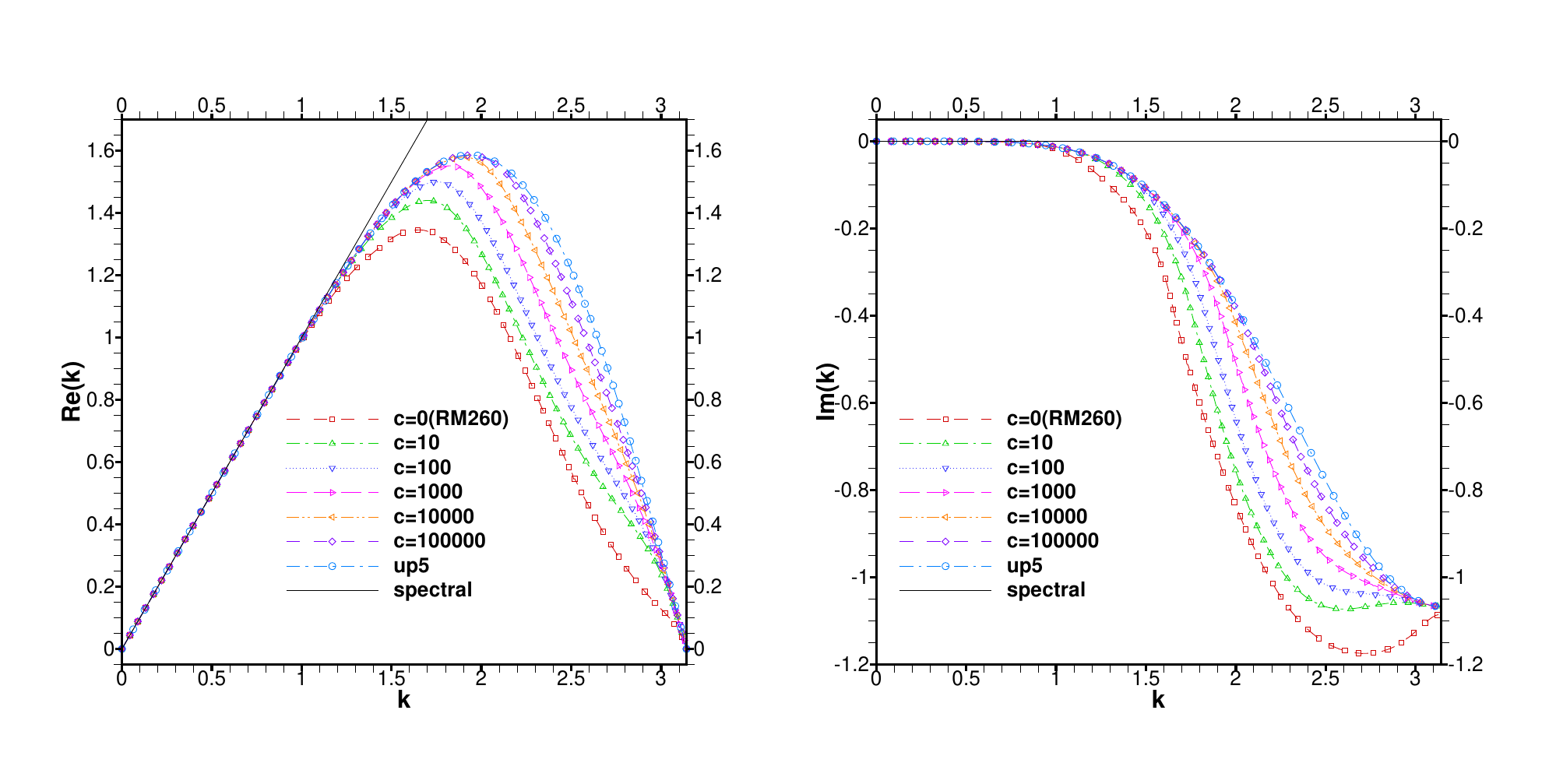}
\small
\caption{Dispersion (left) and dissipation (right) of the WENO-ARMA scheme with selected $c$ when $\chi=100$ by the ADR analysis.}
\end{figure}
\begin{figure}[thb!]
\setlength{\abovecaptionskip}{0pt}
\setlength{\belowcaptionskip}{0pt}
\renewcommand*{\figurename}{Fig.}
\centering
\includegraphics[scale = 0.66]{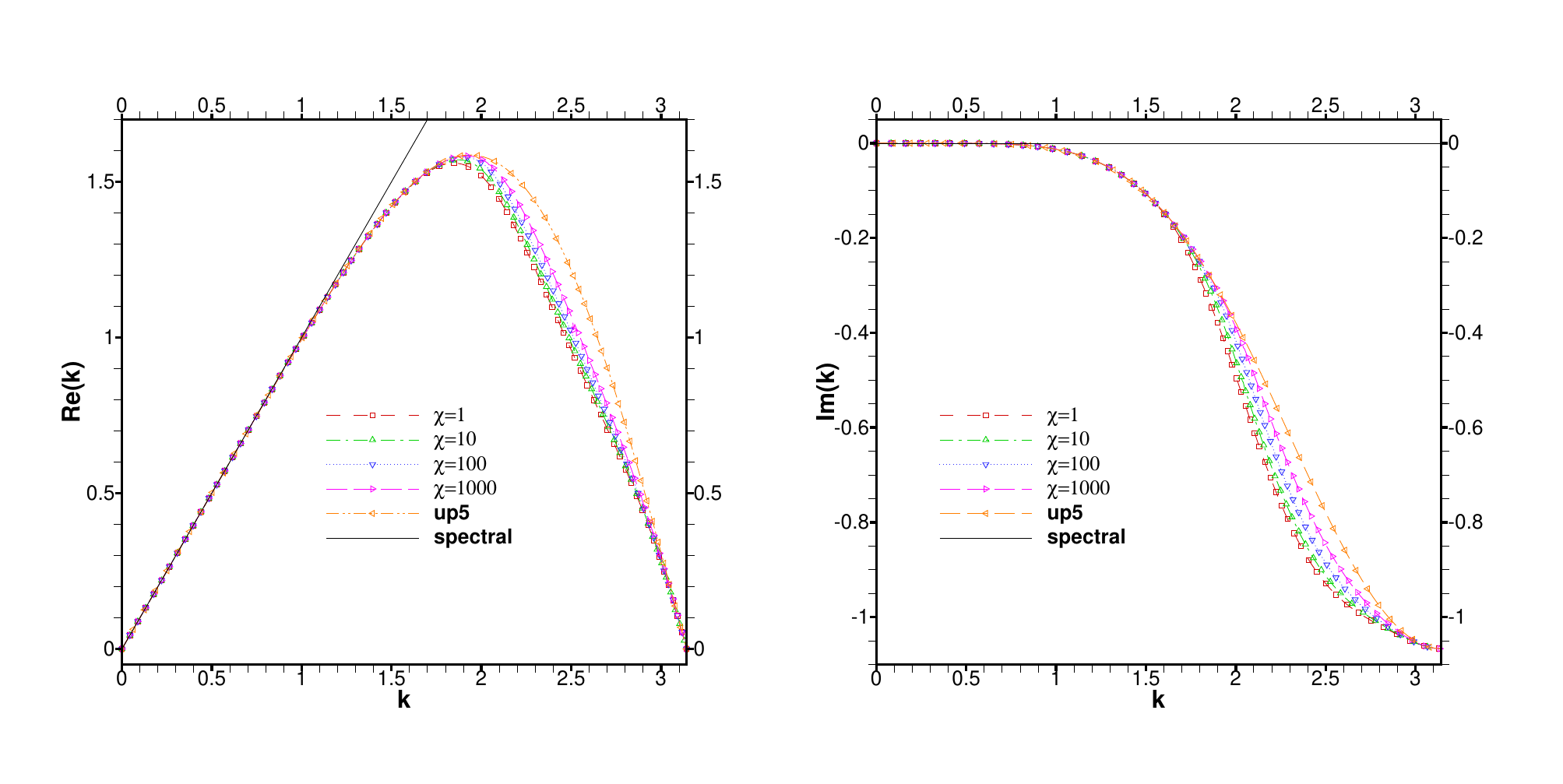}
\small
\caption{Dispersion (left) and dissipation (right) of the WENO-ARMA scheme with selected $\chi$ when $c=10000$ by the ADR analysis.}
\end{figure}
\par
Results of WENO-APMS with selected values of $c$ and $\chi$ by using ADR analysis respectively depicted in Figs. 13 and 14. One can view on these figures that the spectral resolution of WENO-APMS increases with the increase of $c$ or $\chi$. The results of WENO-APMA are like WENO-APMS, and the results are presented in Figs. 15 and 16.
\par
The consequences of WENO-ARMS with selected values of $c$ and $\chi$ by using ADR analysis depicted in Figs. 17 and 18. One can see that the spectral resolution of WENO-ARMS increases as $c$ or $\chi$ increases. Results of WENO-APMA presented in Figs. 19 and 20, and the spectral resolution increase with the increase of $c$ or $\chi$.
\begin{figure}[thb!]
\setlength{\abovecaptionskip}{0pt}
\setlength{\belowcaptionskip}{0pt}
\renewcommand*{\figurename}{Fig.}
\centering
\includegraphics[scale = 0.66]{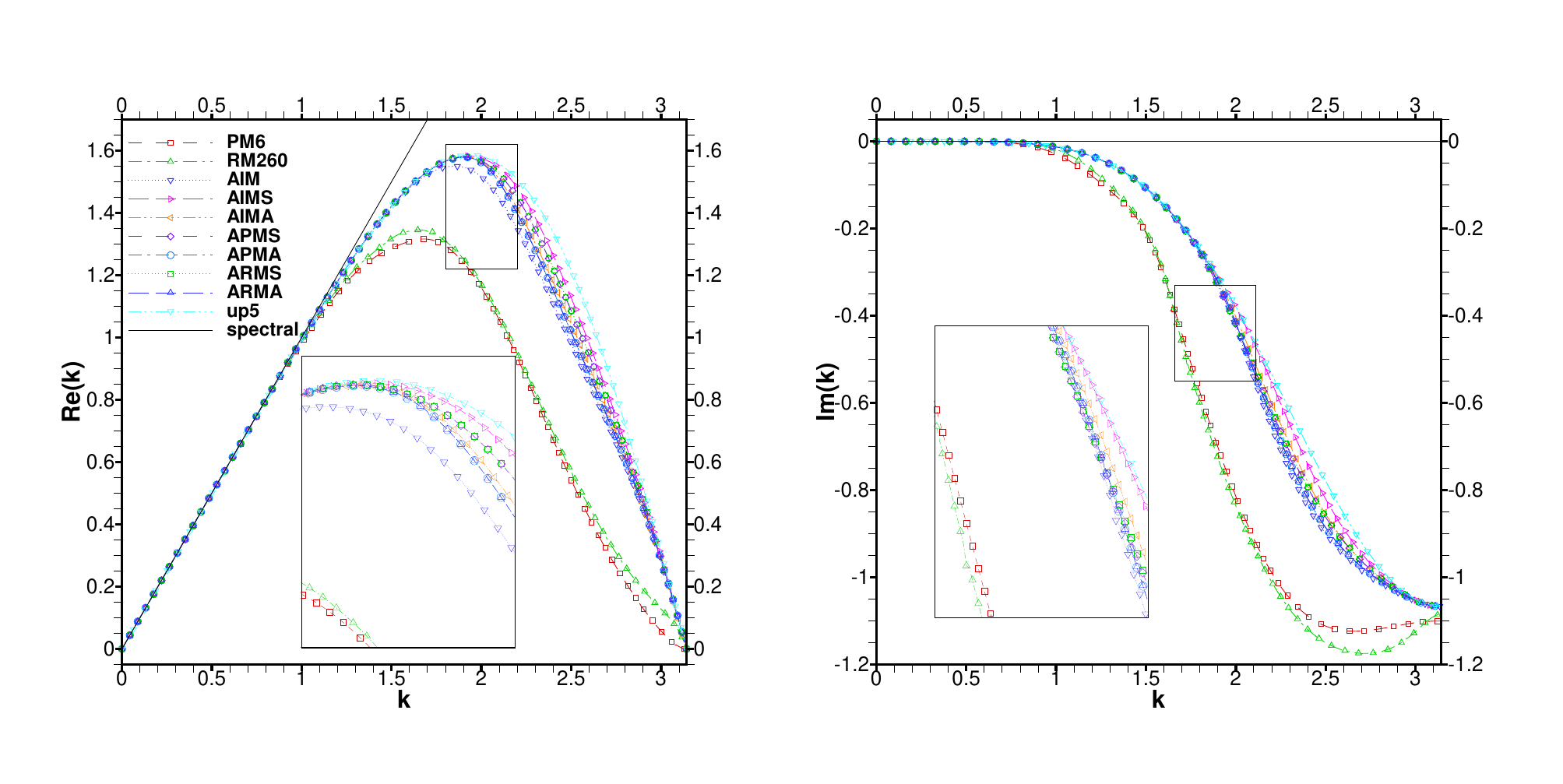}
\small
\caption{Dispersion (left) and dissipation (right) of different WENO schemes by the ADR analysis.}
\end{figure}
\par
Figure 21 shows the comparison of the results of various mapped WENO schemes, where $c$ and $\chi$ in the new mapping functions are taken as $10000$ and $100$, respectively. One can see from the graph, the spectral resolution of all new mapped WENO schemes constructed in this paper is better than WENO-PM6, WENO-RM260, and WENO-AIM. Among them, WENO-AIMS performed best.
\subsection{Determination of the coefficient $\chi$}
\par
From the spectrum analysis results in the previous subsection, we can see that with the further increase of $c$ and $\chi$, all new mapped WENO schemes will eventually evolve into the fifth-order upwind scheme. Considering that when $\chi=0$, both WENO-AIMS, and WENO-AIMA will degenerate into WENO-AIM. For convenience, we only need to give the value of $\chi$ to all new schemes, $c$ and $\kappa$ always taken as $10000$ and $2$. In addition, we only examined the influence of $\chi$ on APMS and APMA with $m=2,~n=6,~\tau=0$, and on APMS and APMA with $n=6$. Therefore, we will determine the value of $\chi$ by solving the one-dimensional wave equation with discontinuities.
\begin{equation}
u_{t}+u_{x}=0,~~x\in[a,~b],
\end{equation}
with periodic boundary conditions, and the initial conditions are set to be
\begin{equation}
u_{0}(x)=\begin{cases}
\frac{1}{6}(G(x,\beta,z-\delta)+G(x,\beta,z+\delta)+4G(x,\beta,z)),&x\in[-0.8,-0.6),\\
1,&x\in[-0.4,-0.2),\\
1-|10(x-0.1)|,&x\in[0,0.2),\\
\frac{1}{6}(F(x,\alpha,a-\delta)+F(x,\alpha,a+\delta)+4F(x,\alpha,z)),&x\in[0.4,0.6),\\
0,&otherwise,
\end{cases}
\end{equation}
where $G(x,\beta,z)=e^{-\beta(x-z)^{2}}$, $F(x,\alpha,a)=\sqrt{max(1-\alpha^{2}(x-a)^{2},0)}$, the constants $z=-0.7$, $\delta=0.0005$, $\beta=ln2/(36\delta^{2})$ and $\alpha=10$.
\par
It is solved on the domain $[-1,1]$ at $t=2,~20,~200$ and $2000~s$. The CFL number was fixed at $0.5$, and $\epsilon=10^{-40}$. The error in the $L^{1}$-norm is comparing the numerical solution $u$ with the exact solution $u_{exact}$, according to
\begin{equation}
\parallel error\parallel_{1}=\frac{1}{N}\sum_{j=0}^{N-1}|u_{j}-(u_{exact})_{j}|
\end{equation}
where $N$ is the number of intervals in space.
\begin{table}
\setlength{\abovecaptionskip}{0pt}
\setlength{\belowcaptionskip}{0pt}\
\footnotesize
\caption{~Performance of WENO-AIMS with (45) at different calculation time for 200 cells.}
\vspace{1 mm}
\begin{center}
\begin{tabular}{c c c c c}
\hline
$\chi$&$t=2~s$&$t=20~s$&$t=200~s$&$t=2000~s$\\
\hline
0(AIM)&$0.23040E-01$&$0.37194E-01$&$0.67009E-01$&$0.11293E+00$\\
1&$0.23011E-01$&$0.37151E-01$&$0.66876E-01$&$0.11314E+00$\\
10&$0.22884E-01$&$0.36969E-01$&$0.66235E-01$&$0.10981E+00$\\
100&$0.22663E-01$&$0.36255E-01$&$0.65301E-01$&$0.10974E+00$\\
1000&$0.22446E-01$&$0.35187E-01$&$0.64512E-01$&$0.11056E+00$\\
\hline
\end{tabular}
\end{center}
\end{table}

\begin{table}
\setlength{\abovecaptionskip}{0pt}
\setlength{\belowcaptionskip}{0pt}\
\footnotesize
\caption{~Performance of WENO-AIMA with (45) at different calculation time for 200 cells.}
\vspace{1 mm}
\begin{center}
\begin{tabular}{c c c c c}
\hline
$\chi$&$t=2~s$&$t=20~s$&$t=200~s$&$t=2000~s$\\
\hline
0(AIM)&$0.23040E-01$&$0.37194E-01$&$0.67009E-01$&$0.11293E+00$\\
1&$0.23006E-01$&$0.37163E-01$&$0.67016E-01$&$0.11021E+00$\\
10&$0.22928E-01$&$0.37139E-01$&$0.67538E-01$&$0.11002E+00$\\
100&$0.22803E-01$&$0.36981E-01$&$0.66891E-01$&$0.10987E+00$\\
1000&$0.22670E-01$&$0.37348E-01$&$0.66771E-01$&$0.11237E+00$\\
\hline
\end{tabular}
\end{center}
\end{table}

\begin{table}
\setlength{\abovecaptionskip}{0pt}
\setlength{\belowcaptionskip}{0pt}\
\footnotesize
\caption{~Performance of WENO-ARMS with (45) at different calculation time for 200 cells.}
\vspace{1 mm}
\begin{center}
\begin{tabular}{c c c c c}
\hline
$\chi$&$t=2~s$&$t=20~s$&$t=200~s$&$t=2000~s$\\
\hline
RM260&$0.24690E-01$&$0.40509E-01$&$0.75446E-01$&$0.11002E+00$\\
1&$0.23198E-01$&$0.36802E-01$&$0.66883E-01$&$0.11216E+00$\\
10&$0.23028E-01$&$0.36189E-01$&$0.62939E-01$&$0.12607E+00$\\
100&$0.22773E-01$&$0.35269E-01$&$0.74297E-01$&$0.13874E+00$\\
1000&$0.22543E-01$&$0.35290E-01$&$0.64657E-01$&$0.13227E+00$\\
\hline
\end{tabular}
\end{center}
\end{table}

\begin{table}
\setlength{\abovecaptionskip}{0pt}
\setlength{\belowcaptionskip}{0pt}\
\footnotesize
\caption{~Performance of WENO-ARMA with (45) at different calculation time for 200 cells.}
\vspace{1 mm}
\begin{center}
\begin{tabular}{c c c c c}
\hline
$\chi$&$t=2~s$&$t=20~s$&$t=200~s$&$t=2000~s$\\
\hline
RM260&$0.24690E-01$&$0.40509E-01$&$0.75446E-01$&$0.11002E+00$\\
1&$0.23158E-01$&$0.37081E-01$&$0.67290E-01$&$0.11620E+00$\\
10&$0.23053E-01$&$0.37192E-01$&$0.67535E-01$&$0.10969E+00$\\
100&$0.22905E-01$&$0.37006E-01$&$0.66053E-01$&$0.11046E+00$\\
1000&$0.22732E-01$&$0.36860E-01$&$0.66544E-01$&$0.11177E+00$\\
\hline
\end{tabular}
\end{center}
\end{table}

\begin{table}
\setlength{\abovecaptionskip}{0pt}
\setlength{\belowcaptionskip}{0pt}\
\footnotesize
\caption{~Performance of WENO-APMS with (45) at different calculation time for 200 cells.}
\vspace{1 mm}
\begin{center}
\begin{tabular}{c c c c c}
\hline
$\chi$&$t=2~s$&$t=20~s$&$t=200~s$&$t=2000~s$\\
\hline
PM6&$0.24795E-01$&$0.40477E-01$&$0.75637E-01$&$0.11002E+00$\\
1&$0.23203E-01$&$0.36656E-01$&$0.66784E-01$&$0.11115E+00$\\
10&$0.23054E-01$&$0.35574E-01$&$0.69388E-01$&$0.13616E+00$\\
100&$0.22784E-01$&$0.35378E-01$&$0.65225E-01$&$0.13189E+00$\\
1000&$0.22546E-01$&$0.35326E-01$&$0.64820E-01$&$0.13227E+00$\\
\hline
\end{tabular}
\end{center}
\end{table}
\begin{table}
\setlength{\abovecaptionskip}{0pt}
\setlength{\belowcaptionskip}{0pt}\
\footnotesize
\caption{~Performance of WENO-APMA with (45) at different calculation time for 200 cells.}
\vspace{1 mm}
\begin{center}
\begin{tabular}{c c c c c}
\hline
$\chi$&$t=2~s$&$t=20~s$&$t=200~s$&$t=2000~s$\\
\hline
PM6&$0.24795E-01$&$0.40477E-01$&$0.75637E-01$&$0.11002E+00$\\
1&$0.23186E-01$&$0.37040E-01$&$0.67248E-01$&$0.11730E+00$\\
10&$0.23064E-01$&$0.37192E-01$&$0.67545E-01$&$0.10969E+00$\\
100&$0.22901E-01$&$0.37007E-01$&$0.66080E-01$&$0.11114E+00$\\
1000&$0.22735E-01$&$0.36859E-01$&$0.66543E-01$&$0.11327E+00$\\
\hline
\end{tabular}
\end{center}
\end{table}
\par
The results for (44) with conditions (45) by using fifth-order WENO-AIMS, AIMA, ARMS, ARMA, APMS, and APMA are shown in Tabs. 1-6. From the results in Table 1, one can see that the error of WENO-AIMS gradually decreases with the increase of $\chi$ as $t\leq 200s$. And $t=2000s$, the error is the smallest as $\chi=100$. Besides, the error of WENO-AIMA that uses asymmetric local operators at $t=2~s$ decreases with an increase of $\chi$, and the smallest error occurs with $\chi=100$ when $t\geq 20~s$. For convenience, in this paper we will take $c=10000$ and $\chi=100$.
\section{Numerical results}
In this section, performances of all adaptive WENO schemes are compared with WENO-AIM, WENO-RM260 and WENO-PM6, and we present the results of our numerical experiments with third-order TVD Runge-Kutta method for time marching\cite{3}
\begin{equation}
\begin{aligned}[b]
u^{(1)}&=u^{n}+\bigtriangleup t^{n}C[u^{n}],\\
u^{(2)}&=\frac{3}{4}u^{n}+\frac{1}{4}(u^{(1)}+\bigtriangleup t^{n}C[u^{(1)}]),\\
u^{(3)}&=\frac{1}{3}u^{n}+\frac{2}{3}(u^{(2)}+\bigtriangleup t^{n}C[u^{(2)}]),
  \end{aligned}
\end{equation}
where $C[u^{n}]$ denote the numerical flux.  \par
In all numerical examples below, the coefficients of adaptive mapping function WENO schemes are selected as: \par
1. WENO-AIMS and WENO-AIMA: $c=1E4,~\kappa=2,~n=4$ and $\chi=100$.\par
2. WENO-ARMS and WENO-ARMA: $c=1E4,~m=2,~n=6,~\tau=0,~\kappa=2$ and $\chi=100$.\par
3. WENO-APMS and WENO-APMA: $c=1E4,~n=6,~\kappa=2$ and $\chi=100$.
\subsection{One dimensional linear advection problems}
\par
The one dimensional linear advection equation (Eq.40) was solved five different initial conditions.
\subsubsection{Linear problem to test the accuracy}
\par
The initial conditions at interval $[-1,~1]$ given as
\begin{equation}
\begin{aligned}[b]
case~1.~u_{0}(x) &= sin(\pi x).\\
case~2.~u_{0}(x) &= sin(\pi x-sin(\pi x)/\pi).
\end{aligned}
\end{equation}
These solutions with critical points are usually used to test the convergence accuracy and numerical dissipation of different schemes. To test the accuracy of various WENO schemes, the calculation time of this problem is assumed to be $t=2~s$. The $L^{1}$-errors, orders and CPU times of WENO scheme with different mapping functions are compared in Table 7. The time steps $\Delta t=(\Delta x)^{5/4}$ is chosen.
\begin{table}
\setlength{\abovecaptionskip}{0pt}
\setlength{\belowcaptionskip}{0pt}\
\footnotesize
\caption{~Comparison of $L^{1}$ errors and orders for Linear problem with case 1 at $t=2~s$.}
\vspace{1 mm}
\begin{center}
\begin{tabular}{c c c c c c c}
\hline
N &WENO-AIM&~~~&WENO-AIMS&~~&WENO-AIMA\\\cline{2-3}\cline{4-5}\cline{6-7}
&Error(order)&~~CPU~time
&Error(order)&~~CPU~time
&Error(order)&~~CPU~time\\
\hline
50&$0.21152E-05(-)$&0.07812&$0.21152E-05(-)$&0.09375&$0.21152E-05(-)$&0.07812\\
100&$0.66058E-07(5.0)$&0.35937&$0.66058E-07(5.0)$&0.37500&$0.66058E-07(5.0)$&0.35937\\
200&$0.20649E-08(5.0)$&1.60938&$0.20649E-08(5.0)$&1.79688&$0.20649E-08(5.0)$&1.62500\\
400&$0.64537E-10(5.0)$&7.48438&$0.64537E-10(5.0)$&8.04688&$0.64537E-10(5.0)$&7.76562\\
800&$0.20170E-11(5.0)$&34.06250&$0.20170E-11(5.0)$&37.31250&$0.20170E-11(5.0)$&36.03125\\
\hline
N &WENO-RM260&~~&WENO-ARMS&~~&WENO-ARMA\\\cline{2-3}\cline{4-5}\cline{6-7}
&Error(order)&~~CPU~time
&Error(order)&~~CPU~time
&Error(order)&~~CPU~time\\
\hline
50&$0.64580E-05(-)$&0.09375&$0.64580E-05(-)$&0.14062&$0.64580E-05(-)$&0.09375\\
100&$0.20168E-06(5.0)$&0.39062&$0.20168E-06(5.0)$&0.48438&$0.20168E-06(5.0)$&0.48438\\
200&$0.63015E-08(5.0)$&1.76562&$0.63015E-08(5.0)$&2.04688&$0.63015E-08(5.0)$&2.03125\\
400&$0.19694E-09(5.0)$&8.48438&$0.19694E-09(5.0)$&9.50000&$0.19694E-09(5.0)$&9.70312\\
800&$0.61552E-11(5.0)$&39.87500&$0.61552E-11(5.0)$&42.93750&$0.61552E-11(5.0)$&44.84375\\
\hline
N &WENO-PM6&~~&WENO-APMS&~~&WENO-APMA\\\cline{2-3}\cline{4-5}\cline{6-7}
&Error(order)&~~CPU~time
&Error(order)&~~CPU~time
&Error(order)&~~CPU~time\\
\hline
50&$0.64580E-05(-)$&0.03125&$0.64580E-05(-)$&0.07812&$0.64580E-05(-)$&0.06250\\
100&$0.20168E-06(5.0)$&0.26562&$0.20168E-06(5.0)$&0.31250&$0.20168E-06(5.0)$&0.31250\\
200&$0.63015E-08(5.0)$&1.14062&$0.63015E-08(5.0)$&1.40625&$0.63015E-08(5.0)$&1.37500\\
400&$0.19694E-09(5.0)$&5.15625&$0.19694E-09(5.0)$&6.29688&$0.19694E-09(5.0)$&6.00000\\
800&$0.61552E-11(5.0)$&24.95312&$0.61552E-11(5.0)$&28.98438&$0.61552E-11(5.0)$&28.28125\\
\hline
\end{tabular}
\end{center}
\end{table}

\begin{table}
\setlength{\abovecaptionskip}{0pt}
\setlength{\belowcaptionskip}{0pt}\
\footnotesize
\caption{~Comparison of $L^{1}$ errors and orders for Linear problem with case 2 at $t=2$.}
\vspace{1 mm}
\begin{center}
\begin{tabular}{c c c c c c c c c}
\hline
N &WENO-AIM&~~~&WENO-AIMS&~~&WENO-AIMA\\\cline{2-3}\cline{4-5}\cline{6-7}
&Error(order)&~~CPU~time
&Error(order)&~~CPU~time
&Error(order)&~~CPU~time\\
\hline
50&$0.24445E-04(-)$&0.09375&$0.24445E-04(-)$&0.10937&$0.24445E-04(-)$&0.10938\\
100&$0.76603E-06(5.0)$&0.39062&$0.76603E-065.0)$&0.39062&$0.76603E-06(5.0)$&0.40625\\
200&$0.23959E-07(5.0)$&1.62500&$0.23959E-07(5.0)$&1.81250&$0.23959E-07(5.0)$&1.71875\\
400&$0.74894E-09(5.0)$&7.68750&$0.74894E-09(5.0)$&8.25000&$0.74894E-09(5.0)$&8.10938\\
800&$0.23410E-10(5.0)$&36.35938&$0.23410E-10(5.0)$&39.37500&$0.23410E-10(5.0)$&38.10938\\
\hline
N &WENO-RM260&~~&WENO-ARMS&~~&WENO-ARMA\\\cline{2-3}\cline{4-5}\cline{6-7}
&Error(order)&~~CPU~time
&Error(order)&~~CPU~time
&Error(order)&~~CPU~time\\
\hline
50&$0.74488E-04(-)$&0.10937&$0.74498E-04(-)$&0.06250&$0.74498E-04(-)$&0.10937\\
100&$0.23362E-05(5.0)$&0.37500&$0.23362E-05(5.0)$&0.43750&$0.23362E-05(5.0)$&0.37500\\
200&$0.73077E-07(5.0)$&1.76562&$0.73077E-07(5.0)$&1.90625&$0.73077E-07(5.0)$&1.89062\\
400&$0.22856E-08(5.0)$&8.21875&$0.22856E-08(5.0)$&9.40625&$0.22856E-08(5.0)$&9.26562\\
800&$0.71437E-10(5.0)$&39.35938&$0.71437E-10(5.0)$&44.14062&$0.71437E-10(5.0)$&43.37500\\
\hline
N &WENO-PM6&~~&WENO-APMS&~~&WENO-APMA\\\cline{2-3}\cline{4-5}\cline{6-7}
&Error(order)&~~CPU~time
&Error(order)&~~CPU~time
&Error(order)&~~CPU~time\\
\hline
50&$0.76265E-04(-)$&0.07812&$0.74498E-04(-)$&0.09375&$0.74498E-04(-)$&0.06250\\
100&$0.23410E-05(5.0)$&0.26562&$0.23362E-05(5.0)$&0.31250&$0.23362E-05(5.0)$&0.34375\\
200&$0.73085E-07(5.0)$&1.21875&$0.73077E-07(5.0)$&1.34375&$0.73077E-07(5.0)$&1.42188\\
400&$0.22856E-08(5.0)$&5.64062&$0.22856E-08(5.0)$&6.65625&$0.22856E-08(5.0)$&6.89062\\
800&$0.71437E-10(5.0)$&26.62500&$0.71437E-10(5.0)$&31.48438&$0.71437E-10(5.0)$&30.07812\\
\hline
\end{tabular}
\end{center}
\end{table}
\par
One can see from the table, all of WENO schemes using different mapping functions have achieved fifth-order accuracy. But we also noticed that the two improved WENO-AIM schemes have the same errors as WENO-AIM when calculating the two problems, and the CPU time is slightly longer than the original formula. Similar results also appear in other improved mapping WENO schemes.
\subsubsection{Linear problem with discontinuities}
The initial conditions are given as
\begin{equation}\begin{aligned}[b]
case3.~u_{0}(x)&=\begin{cases}
1,&x\in[-1,0),\\
0,&x\in[0,1].
\end{cases}\\
case4.~u_{0}(x)&=\begin{cases}
-sin(\pi x)-\frac{1}{2}x^{3},&x\in[-1,0],\\
-sin(\pi x)-\frac{1}{2}x^{3}+1,&x\in(0,1],
\end{cases}\\
case5.~ u_{0}(x) &~defined~by~(45).
\end{aligned}
\end{equation}
\begin{figure}[thb!]
\centering
\setlength{\abovecaptionskip}{0pt}
\setlength{\belowcaptionskip}{1pt}
\renewcommand*{\figurename}{Fig.}
\includegraphics[scale = 0.55]{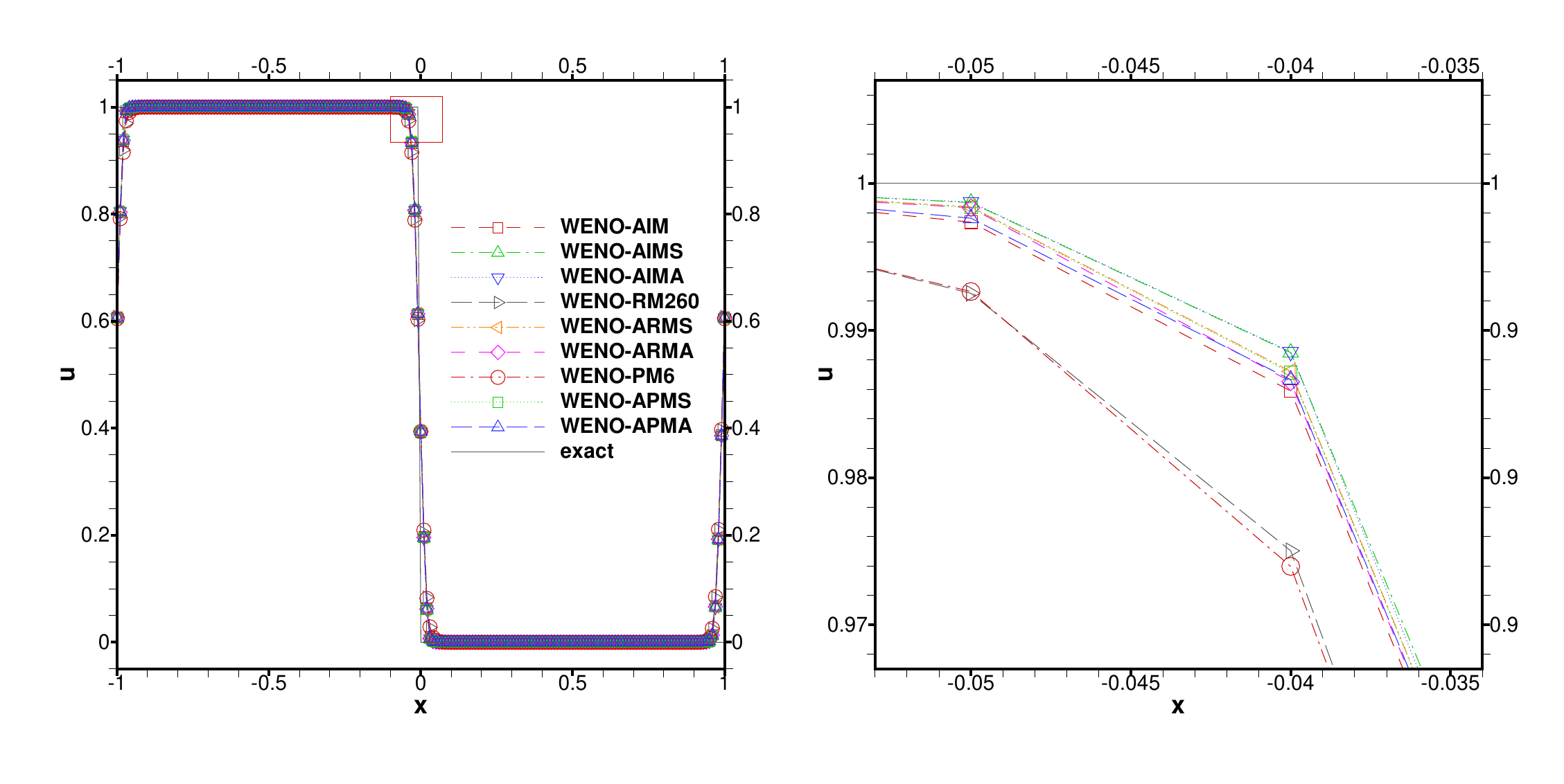}
\small
\caption{Performance of fifth-order mapped WENO schemes for case 3 at $t=2~s$ with 200 grids.}
\end{figure}
\begin{figure}[thb!]
\centering
\setlength{\abovecaptionskip}{0pt}
\setlength{\belowcaptionskip}{1pt}
\renewcommand*{\figurename}{Fig.}
\includegraphics[scale = 0.55]{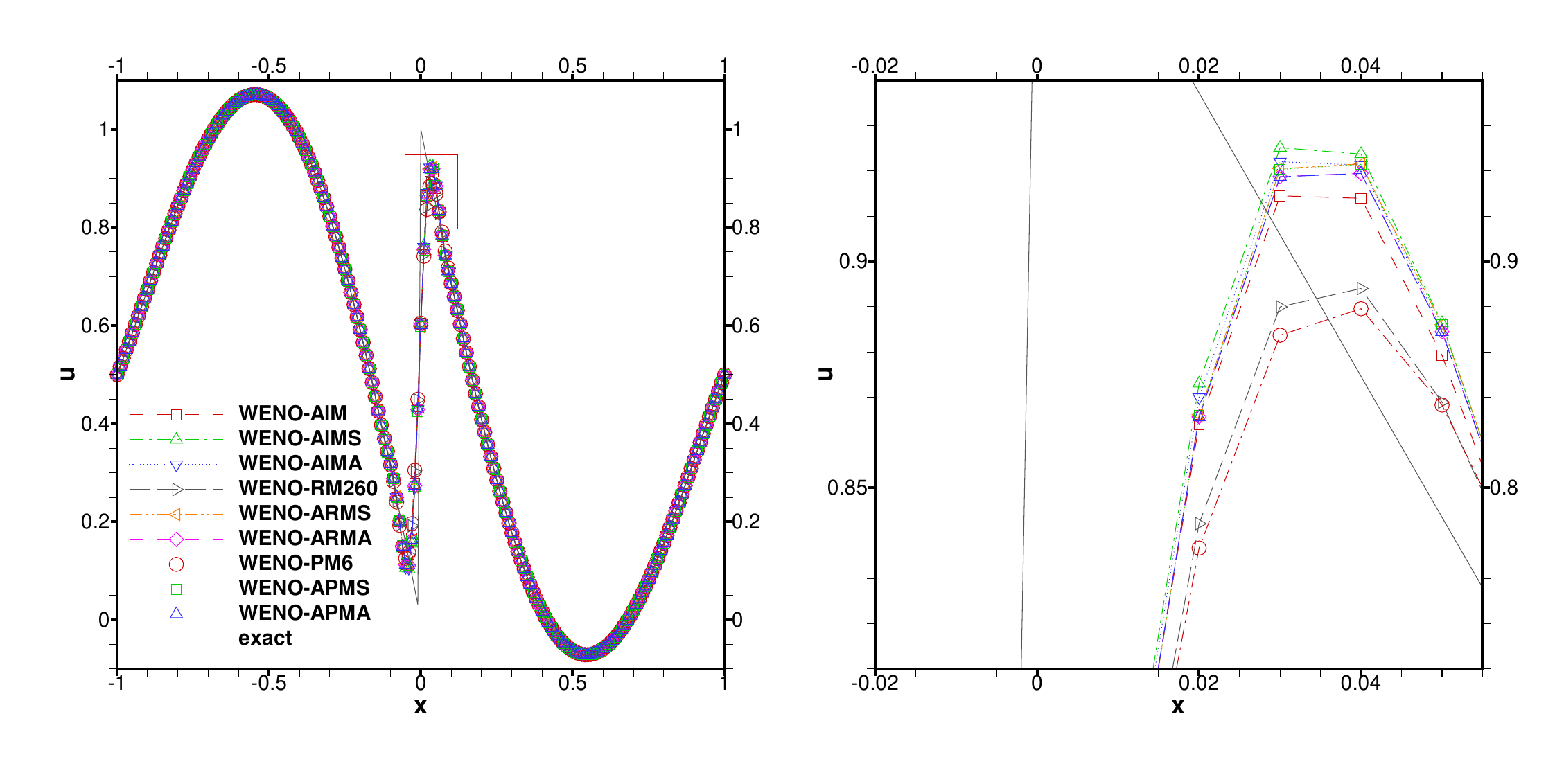}
\small
\caption{Performance of fifth-order mapped WENO schemes for case 4 at $t=2~s$ with 200 grids.}
\end{figure}
\begin{figure}[thb!]
\centering
\setlength{\abovecaptionskip}{0pt}
\setlength{\belowcaptionskip}{1pt}
\renewcommand*{\figurename}{Fig.}
\includegraphics[scale = 0.55]{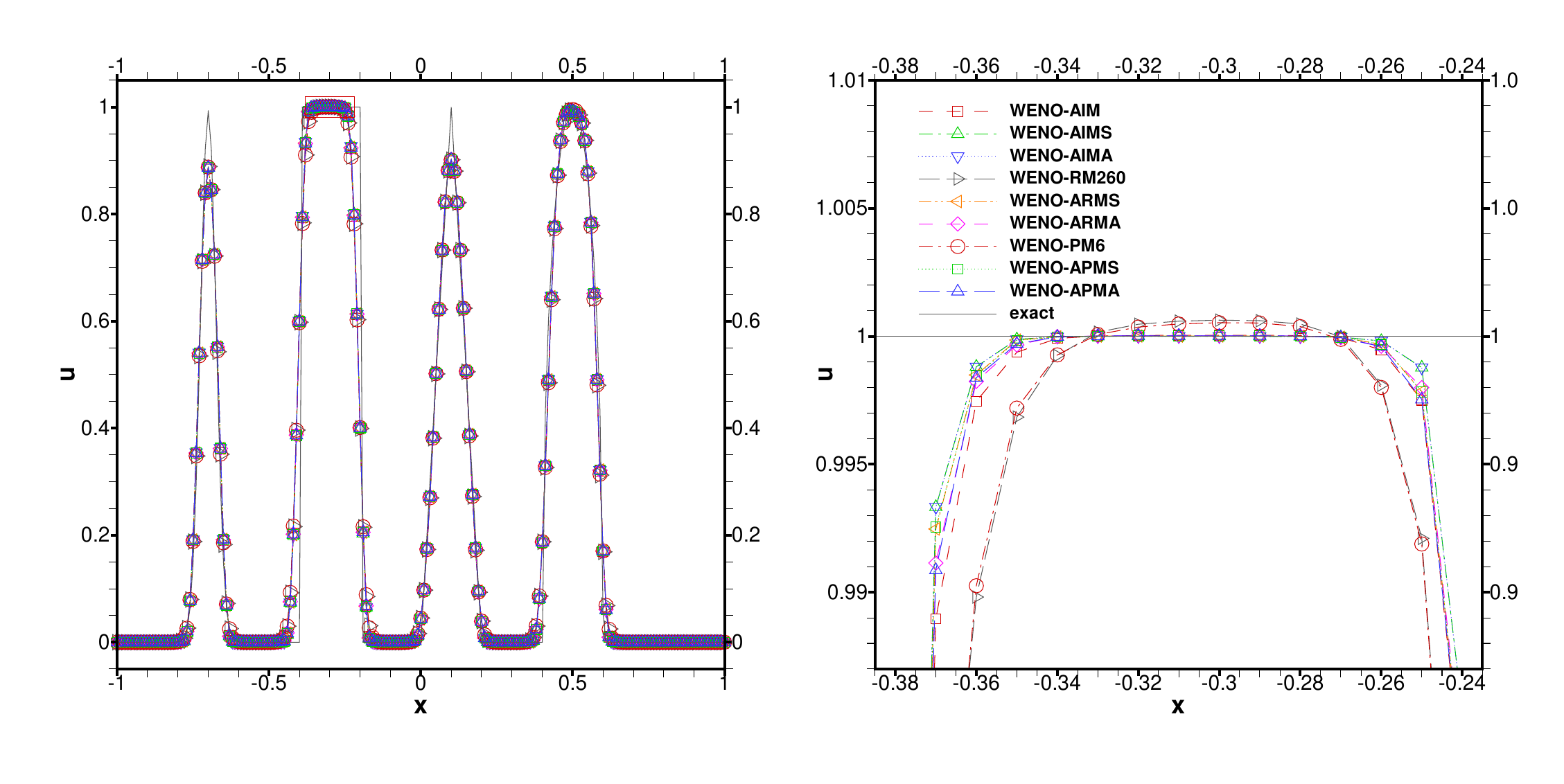}
\small
\caption{Performance of fifth-order mapped WENO schemes for case 5 at $t=2~s$ with 200 grids.}
\end{figure}
\par
Figures. 21-23 depict respectively the numerical solution of all nine WENO schemes at $t=2~s$ with 200 grids by calculating the one-dimensional linear advection equation (Eq.40) with initial conditions (49), (50), and (45). The time steps selected by $\Delta t=0.6(\Delta x)$. As seen from these figures, these mapped WENO schemes proposed in this paper perform better than the original WENO-AIM, WENO-RM, and WENO-PM. And in the six new WENO schemes, both WENO-AIMS and WENO-AIMA are the same and perform better than the other four.
\subsection{One-dimensional Euler problems}
Consider the following nonlinear 1-D Euler system
\begin{gather}\textbf{u}_{t}+f(\textbf{u})_{x}=0,\end{gather}
with
\begin{gather}\textbf{u}=(\rho,\rho u, E)^{T},~~f(\textbf{u})=(\rho u,\rho u^{2}+p, u(E+p))^{T},~E=\frac{p}{\gamma-1}+\frac{1}{2}\rho u^{2},\end{gather}
\noindent
where $\rho,~u,~p$ and $e$ are the density, velocity, pressure and total energy, respectively. Time step is taken as
\begin{equation}
\Delta t = \frac{CFL \Delta x}{max_{i}(u_{i}+c_{i})}
\end{equation}
where $CFL=0.5$, $c$ is the speed of sound and given as $c=\sqrt{\gamma p/\rho}$.
\subsubsection{SOD's problem}
The initial conditions of the problem are as follows\cite{15}
\begin{equation}
(\rho, u, p)=\begin{cases}
(0.125, 0, 0.1), &x\in[0,0.5],\\
(1.000, 0, 1.0), &x\in(0.5,1].
\end{cases}
\end{equation}
The initial conditions for this problem include a fixed gas with a right-traveling shock wave, contact discontinuities, and a left-traveling expansion wave containing sound velocity points, which can usually test the performance of different schemes. Fig. 25 gives the distribution of density at $t = 0.14~s$ on a 200-cell with the reference solution calculated on a 1000-cell grid using the WENO-AIM. One can see from the figure, all solutions can smoothly pass through the sound velocity point, and other wave systems are also close to the reference solution. Among them, all improved WENO schemes perform better than WENO-AIM, WENO-RM260, and WENO-PM6. But we can see that there is an overshoot in the results of WENO-AIMS. In contrast, the result of WENO-AIMA is closer to the exact solution than other improved WENO schemes.
\begin{figure}[thb!]
\centering
\setlength{\abovecaptionskip}{0pt}
\setlength{\belowcaptionskip}{0pt}
\renewcommand*{\figurename}{Fig.}
\includegraphics[scale = 0.525]{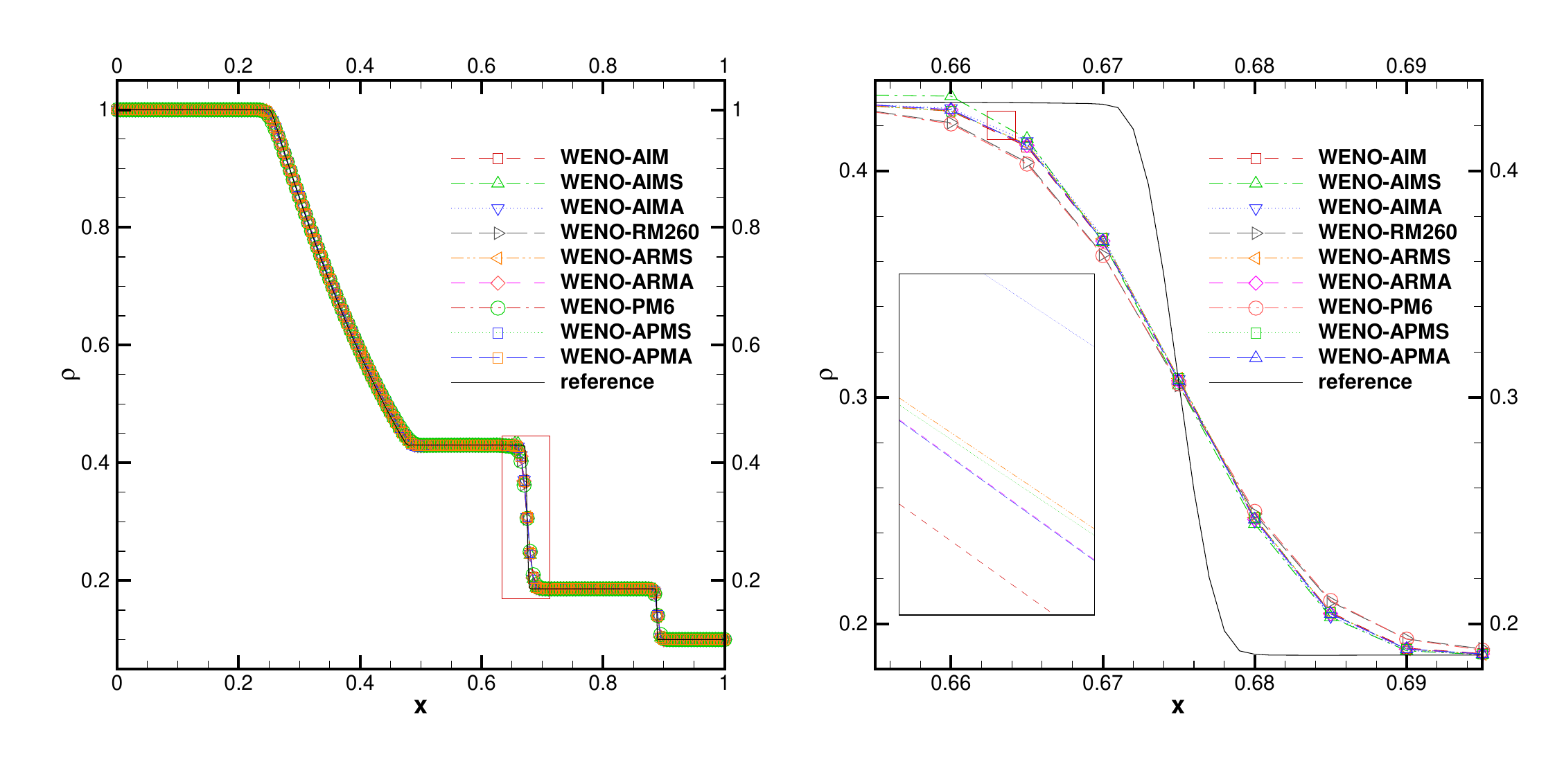}
\small
\caption{Performance of fifth order WENO schemes for SOD's problem at $t=0.14~s$ with 200 grids.}
\end{figure}
\subsubsection{Shu-Osher's problem}
The initial conditions for the problem are as follows\cite{3}
\begin{equation}
(\rho, u, p)=\begin{cases}
(3.857143, 2.629369, 10.33333), &x\in[-5,-4),\\
(0.125, 0, 0.1), &x\in[-4,5].
\end{cases}
\end{equation}
This problem includes low-frequency and high-frequency density disturbances and is often used as a benchmark case for testing the performance of different schemes. Fig.26 presents the distribution of density at $t=1.8~s$ on a uniform 200 cells with the reference solution calculated on a 1000-cell grid using the WENO-IM. As seen from the figure, all the new WENO schemes perform slightly better than WENO-AIM, better than WENO-RM260 and WENO PM6.
\begin{figure}[thb!]
\centering
\setlength{\abovecaptionskip}{0pt}
\setlength{\belowcaptionskip}{0pt}
\renewcommand*{\figurename}{Fig.}
\includegraphics[scale = 0.525]{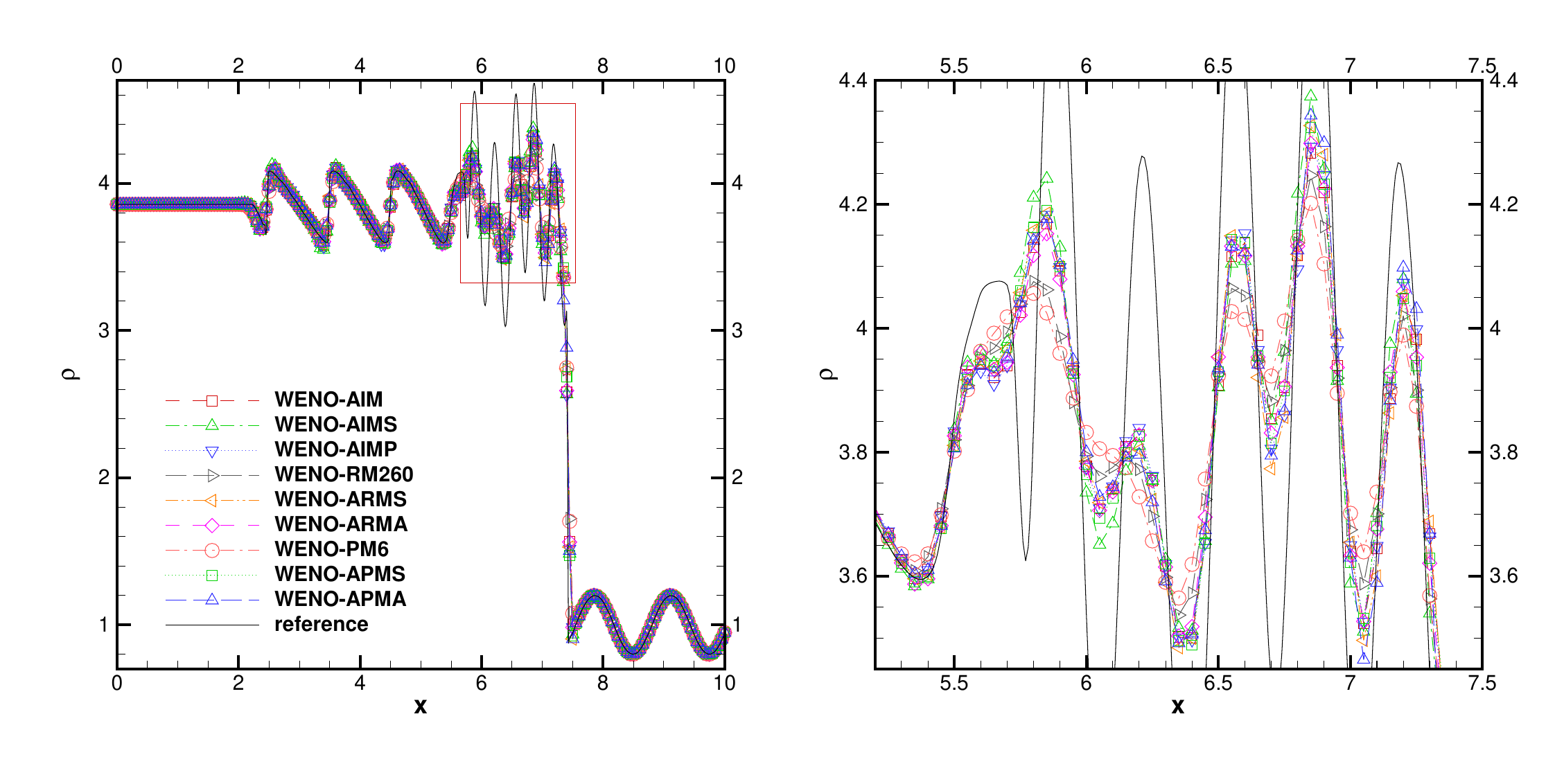}
\small
\caption{Performance of fifth order WENO schemes for Shu-Osher's problem at $t=1.8~s$ with 200 grids.}
\end{figure}
\subsection{Two dimension problem}
In this section, we mainly investigate the performance characteristics of different WENO schemes by solving two-dimensional problems. The 2-D Euler equations are given as
\begin{equation}
\textbf{u}_{t}+f(\textbf{u})_{x}+g(\textbf{u})_{y}=0,
\end{equation}
with
\begin{equation}
\begin{aligned}[b]
\textbf{u}&=(\rho,\rho u, \rho v ,E)^{T},\\
f(\textbf{u})&=(\rho u,\rho u^{2}+p, \rho uv, u(E+p))^{T},\\
g(\textbf{u})&=(\rho v,\rho uv, \rho v^{2}+p, v(E+p))^{T},
  \end{aligned}
\end{equation}
where, $\rho$ is the density, $u$ and $v$ are $x-$ and $y-$velocities, respectively. $E$ is the total energy, $p$ is the pressure, which is related to the energy by $E=\frac{p}{\gamma-1}+\frac{1}{2}\rho(u^{2}+v^{2})$ with the CFL number taken as $0.5$, and the specific heat ratio $\gamma=1.4$ for all examples.
\subsubsection{Periodic vortex propagation Problem}
This is a 2-D periodic problem used to assess the numerical dissipation of the finite difference schemes. The vortex is described as a perturbation to the mean flow, and the initial condition is\cite{18,19}
\begin{equation}
\begin{cases}\rho=\left[1-\frac{(\gamma-1)\beta^{2}}{8\gamma\pi^{2}}e^{1-r^{2}}\right]^{\frac{1}{\gamma-1}},\\
u=1-\frac{\beta}{2\pi}e^{(1-r^{2})/2}\overline{y},\\
v=1+\frac{\beta}{2\pi}e^{(1-r^{2})/2}\overline{x},\\
p=\rho^{\gamma},
\end{cases}
\end{equation}
where $r=\sqrt{\overline{x}^{2}+\overline{y}^{2}}$, $(\overline{x},\overline{y})=(x-5,y-5)$ and $\beta=5$.
\begin{figure}[thb!]
\setlength{\abovecaptionskip}{0pt}
\setlength{\belowcaptionskip}{0pt}
\renewcommand*{\figurename}{Fig.}
\begin{center}
\includegraphics[scale = 0.35]{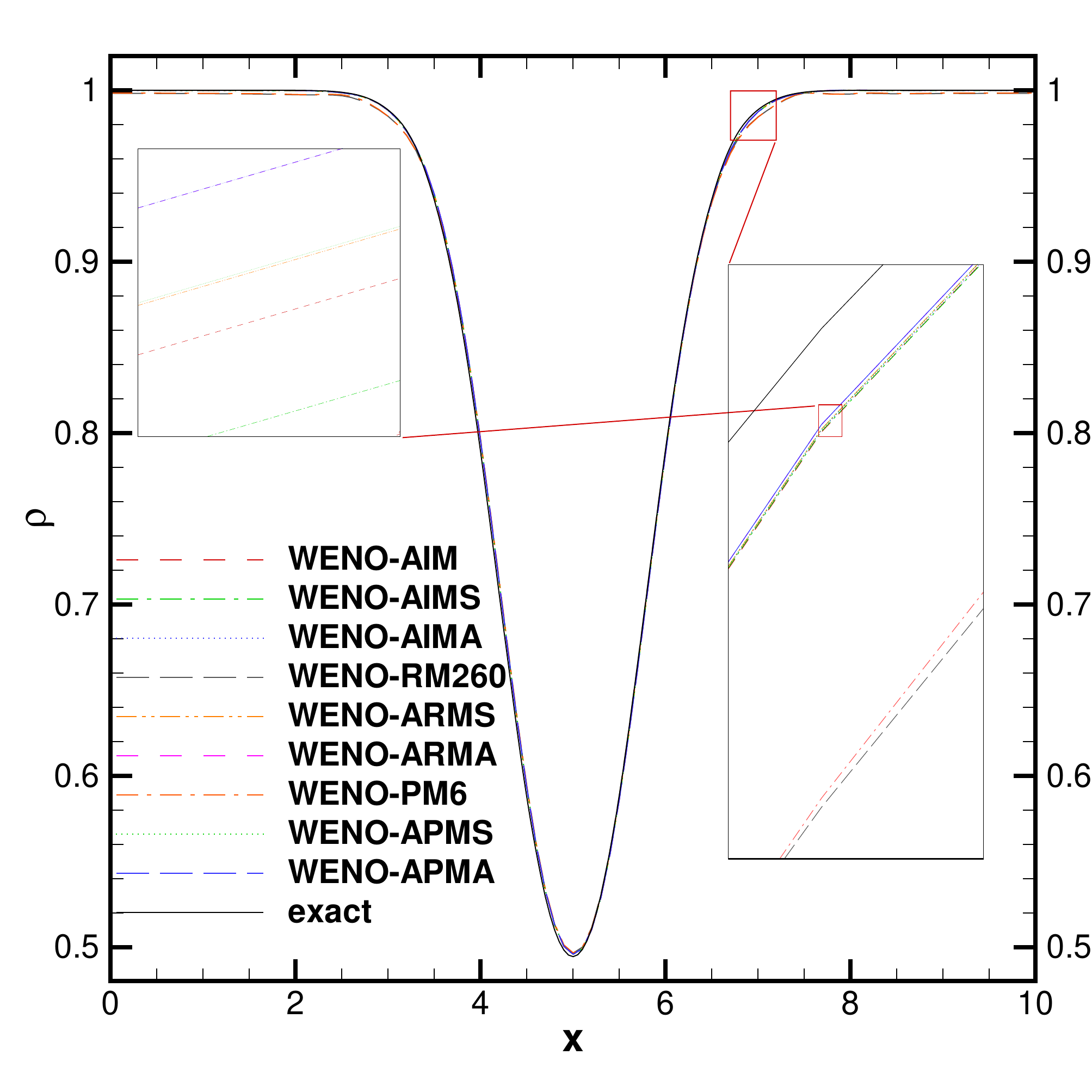}
\end{center}
\small
\caption{The distribution of pressure along $y=5$ using fifth order WENO schemes for periodic vortex propagation problem at $t=100~s$ with $100\times 100$ grids.}
\end{figure}\par
The computational region is $[0, 10] \times [0, 10]$ and periodic boundary condition used at all four boundaries. Fig. 27 presents the solution of nine WENO schemes until the ultimate time $t=100~s$ using uniform meshes with $100\times 100$ grids. One can see from the partially enlarged view of the graph that the effects of several improved schemes are better than those of WENO-AIM, WENO-RM260, and WENO-PM6. Among them, all new mapped WENO schemes equip with asymmetric local operators perform the same and better than others.
\subsubsection{2-D Riemann Problem}
The initial conditions for 2-dimensional Riemann problem are as following\cite{21}
\begin{equation}
(\rho, u, v, p)=\begin{cases}
(1, 0.1, 1, 0.1), &\text{0.6 $\leq x \leq $1, 0.6 $\leq y \leq $1},\\
(0.5313, 0.8276, 0, 0.4), &\text{0 $\leq x < $0.6, 0.6 $\leq y \leq $1},\\
(0.8, 0.1, 0, 0.4), &\text{0 $\leq x < $0.6, 0 $\leq y < $0.6},\\
(0.5313, 0.1, 0.7276, 0.4), &\text{0.6 $< x \leq $1, 0 $\leq y < $0.6}.
\end{cases}
\end{equation}
The calculation region is $[0,~1]\times[0,~1]$, and zero-order extrapolation boundary conditions assumed on all boundaries. Figure 28 depicts the results of nine WENO schemes up to the last time $t = 0.8~s$ using a uniform grid of $400\times 400$ with 12 equidistant density contours ranging from 0.2 to 1.7.
\begin{figure}[thb!]
\centering
\setlength{\abovecaptionskip}{0pt}
\setlength{\belowcaptionskip}{0pt}
\renewcommand*{\figurename}{Fig.}
\includegraphics[scale = 0.7]{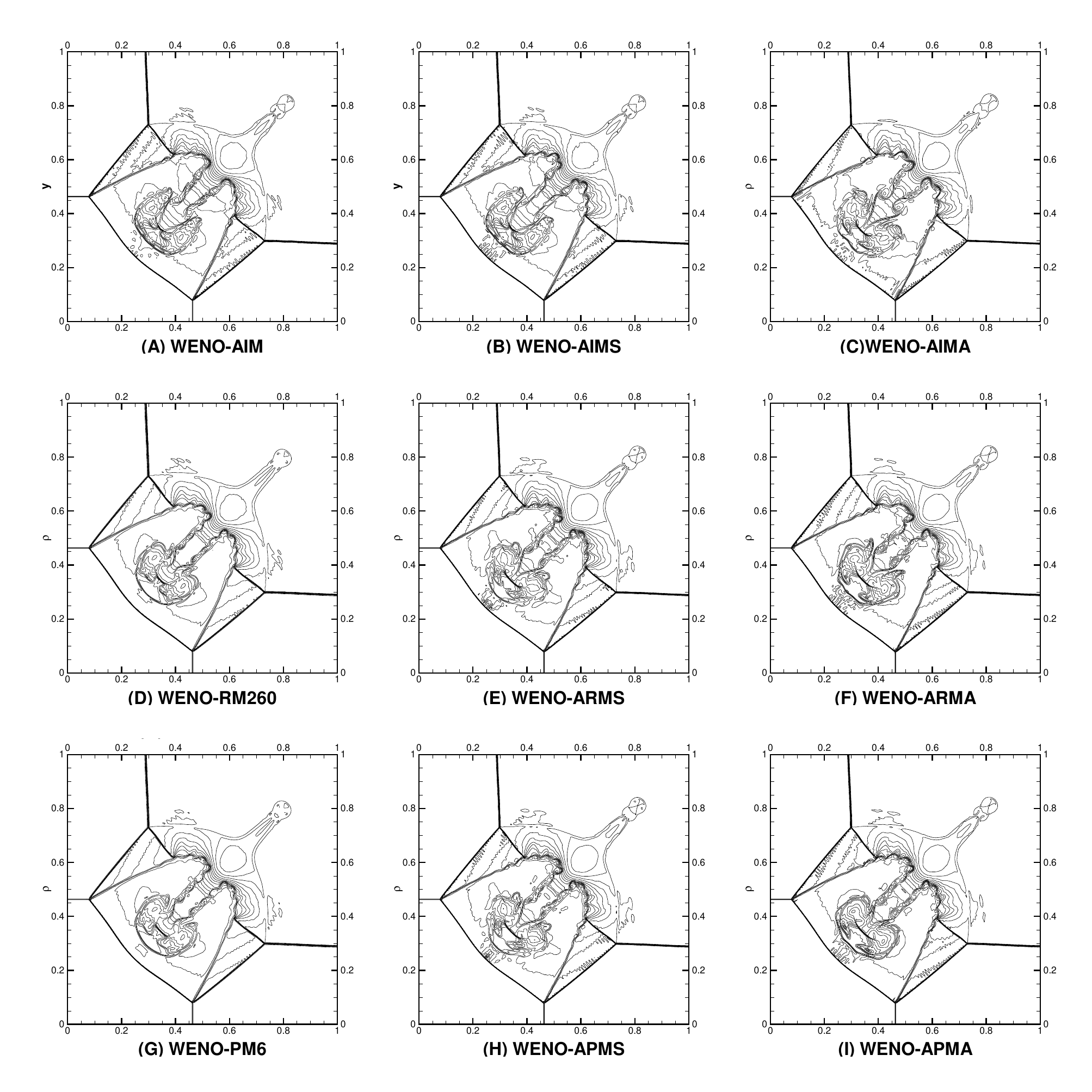}
\small
\caption{Performance of fifth order WENO schemes for 2-D Riemann problem at $t=0.8~s$ with $400\times 400$ grids,~18 contours from 0.2 to 1.7.}
\end{figure}\par
As seen from these figures, all mapped WENO schemes can capture reflected shocks and contact discontinuities satisfactorily. And the WENO-AIMS, WENO-APMS, and WENO-APMA can construct finer unstable Kelvin-Helmholtz microstructure than others. Among these results of all new WENO schemes, WENO-AIMS, WENO-AIMA, WENO-APMS, and WENO APMA perform the same and better than the others. But we also see that at the mushroom-shaped front edge, the three improved WENO schemes using asymmetric local operators appear asymmetry.
\subsubsection{Implosion Problem}
This problem is often used to test the ability of numerical schemes to resolve contact discontinuities and maintain symmetry: if the scheme does not maintain symmetry, no jets will be generated or distorted\cite{21}. The initial conditions
\begin{equation}
(\rho, u, v, p)=\begin{cases}
(0.125, 0, 0, 0.14), &\text{0 $\leq |x|+|y| < $ 0.15},\\
(1, 0,  0, 1), &\text{otherwise}.
\end{cases}
\end{equation}
The computational region is $[-0.3, 0.3] \times [-0.3, 0.3]$, reflecting boundary conditions on all four boundaries. Fig. 29 presents the solution of nine WENO schemes up to $t=2.5~s$ using a uniform grid of  $400\times 400$ with 10 equidistant density contours ranging from 0.45 to 1.05.
\begin{figure}[thb!]
\centering
\setlength{\abovecaptionskip}{0pt}
\setlength{\belowcaptionskip}{0pt}
\renewcommand*{\figurename}{Fig.}
\includegraphics[scale = 0.7]{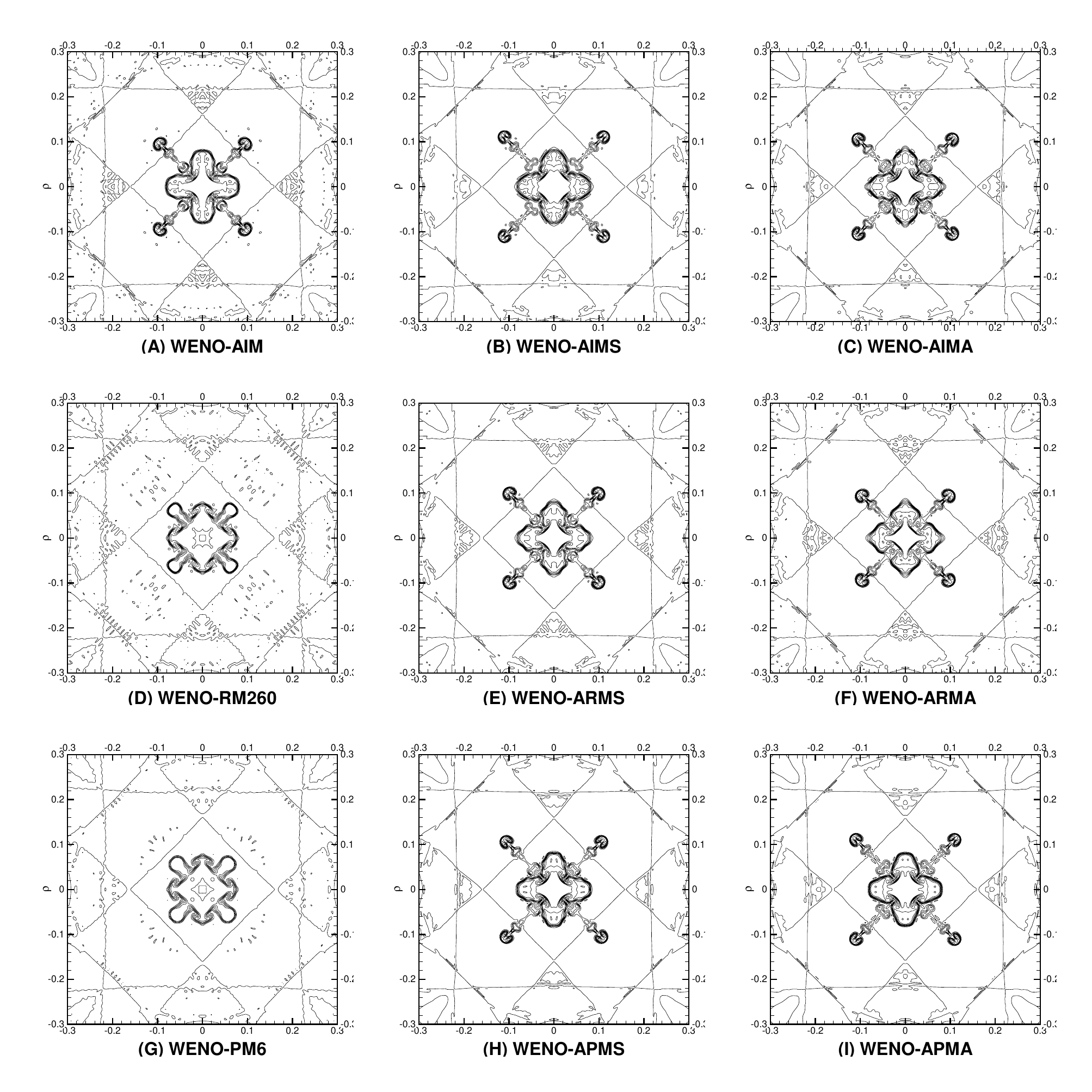}
\small
\caption{Performance of fifth order WENO schemes for implosion problem at $t=2.5~s$ with $400\times 400$ grids,~10 contours from 0.45 to 1.05.}
\end{figure}\par
From the figure, all WENO schemes can capture the jet and maintain symmetry well. The jets of all the improved schemes are longer than those of the corresponding schemes, which shows that the improved schemes have higher resolution. Among all new WENO schemes, WENO-AIMS, WENO-AIMA, WENO-APMS, and WENO APMA perform the same and better than the others.
\subsubsection{Double-Mach reflection problem}
This problem has widely used as a test example for high-order schemes. The initial conditions as following\cite {3,13}
\begin{equation}
(\rho,~u,~v,~p)=\begin{cases}
(8,57.1597,-33.0012,563.544),&y< \sqrt{3}(x-1/6),\\
(1.4,0,0,2.5),&y\geq \sqrt{3}(x-1/6),\\
\end{cases}
\end{equation}
where the computational region is $[0, 4]\times [0, 1]$. For the bottom, the exact post-shock condition is imposed at the interval $[0,0.6]$, and the reflective boundary condition is used for the rest. The exact motion of a Mach 10 shock is set to the top boundary. Inflow and outflow boundary conditions are set for the left and right boundaries, respectively.
\begin{figure}[thb!]
\centering
\setlength{\abovecaptionskip}{0pt}
\setlength{\belowcaptionskip}{0pt}
\renewcommand*{\figurename}{Fig.}
\includegraphics[scale = 0.7]{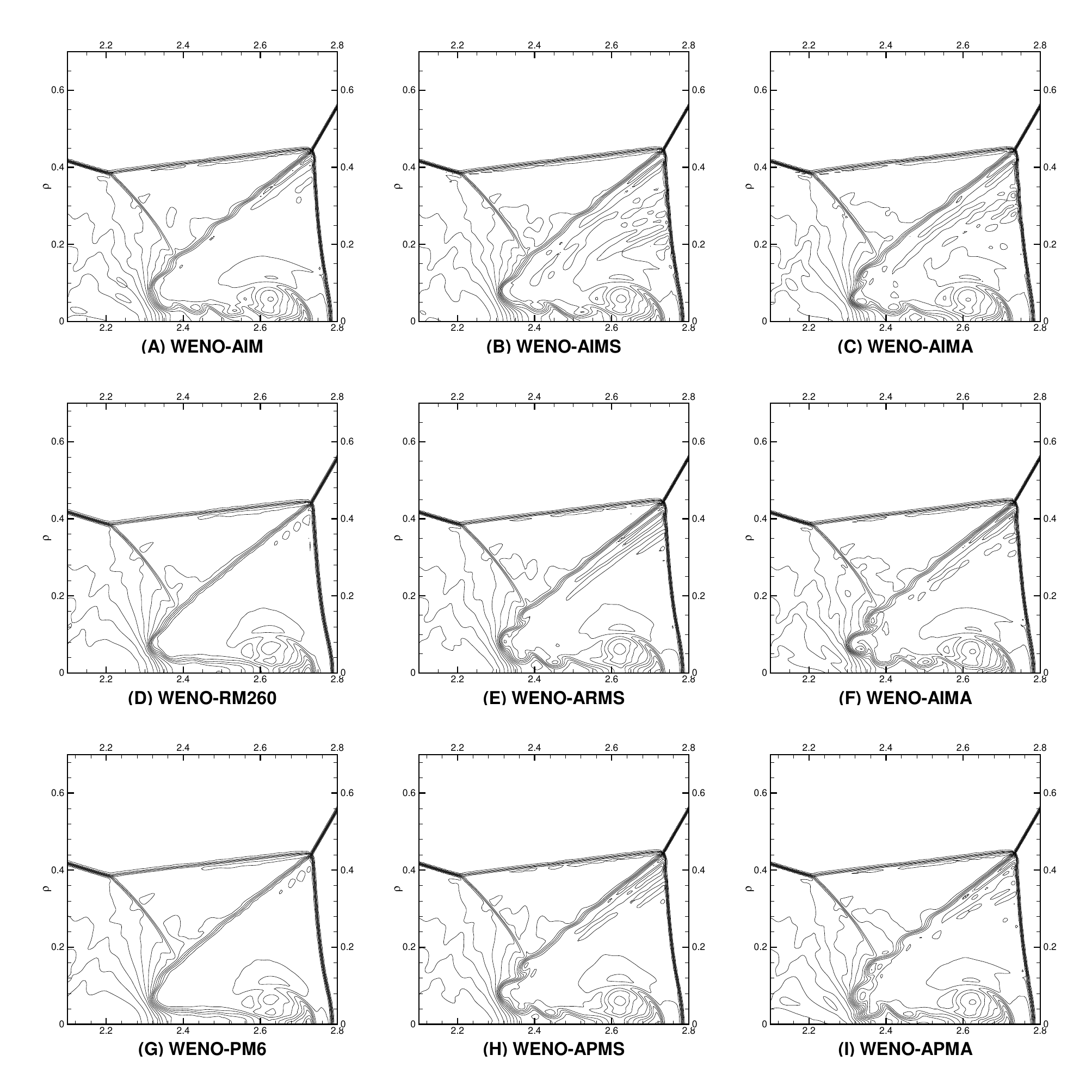}
\small
\caption{Performance of fifth order WENO schemes for double-Mach reflection problem at $t=0.2~s$ with $960\times 240$ grids,~30 contours from 2 to 22.}
\end{figure}
\par
Fig. 30 presents the solution of all mapped WENO schemes up to the final time $t=0.2~s$ using $960\times 240$ uniform meshes with 30 equidistant density contours ranging from 2 to 22. One can see that WENO schemes with all mapping functions can capture shock structures well, and WENO-ARMA and WENO-APMA schemes can resolve more instability microstructure than others.
\subsubsection{computational efficiency}
Table 9 shows the CPU time of each scheme in different 2D situations. We use the same desktop workstation to calculate all the tests. As seen from the results that the computational time of WENO-RM, WENO-ARMS, and WENO-ARMP are much longer than other schemes. And WENO-AIMS, WENO-AIMA, WENO-APMS, and WENO-APMA have higher efficiency.
\begin{table}
\setlength{\abovecaptionskip}{0pt}
\setlength{\belowcaptionskip}{0pt}\
\footnotesize
\caption{~The average computational time of all 2D problems. The normalized results of other schemes relative to the computational time of WENO-AIM are given in parentheses.}
\vspace{1 mm}
\begin{center}
\begin{tabular}{c c c c c}
\hline
Case&Grid number&WENO-AIM&WENO-AIMS&WENO-AIMA\\
\hline
Vortex problem&$100\times100$&$0.6392E+03(1.00)$&$0.6779E+03(1.06)$&$0.6720E+03(1.05)$\\
Riemman problem&$400\times400$&$0.2594E+04(1.00)$&$0.2831E+04(1.09)$&$0.2795E+04(1.07)$\\
Implosion&$400\times400$&$0.1464E+05(1.00)$&$0.1957E+05$(1.33)&$0.1507E+05(1.03)$\\
Double Mach&$960\times240$&$0.3424E+04(1.00)$&$0.3786E+04(1.10)$&$0.3742E+04(1.09)$\\
\hline
\hline
Case&Grid number&WENO-RM260&WENO-ARMS&WENO-ARMA\\
\hline
Vortex problem&$100\times100$&$0.1137E+04(1.78)$&$0.1340E+04(2.09)$&$0.1291E+04(2.01)$\\
Riemman problem&$400\times400$&$0.7164E+04(2.76)$&$0.7978E+04(3.07)$&$0.7520E+04(2.90)$\\
Implosion problem&$400\times400$&$0.3553E+05(2.42)$&$0.3171E+05(2.16)$&$0.3203E+05(2.19)$\\
Double Mach problem&$960\times240$&$0.4968E+04$&$0.5657E+04$&$0.5749E+04$\\
\hline
\hline
Case&Grid number&WENO-PM6&WENO-APMS&WENO-APMA\\
\hline
Vortex problem&$100\times100$&$0.6865E+03(1.07)$&$0.7731E+03(1.21)$&$0.7634E+03(1.20)$\\
Riemman problem&$400\times400$&$0.2859E+04(1.10)$&$0.4025E+04(1.55)$&$0.2969E+04(1.14)$\\
Implosion&$400\times400$&$0.1666E+05(1.14)$&$0.2245E+05(1.53)$&$0.2237E+05(1.54)$\\
Double Mach&$960\times240$&$0.2753E+04(0.80)$&$0.3094E+04(0.90)$&$0.3057E+04(0.89)$\\
\hline
\end{tabular}
\end{center}
\end{table}
\section{Conclusions}
\label{}
\vspace{3 mm}
\par
In this paper, by improving the mapping function of the WENO-AIM scheme with two local operators, we composed two more generalized adaptive mapped WENO schemes, WENO-AIMS and WENO-AIMA. We also introduced these two local operators into the mapping functions of WENO-RM and WENO-PM and composed the adaptive WENO-RM and adaptive WENO-PM schemes. In the numerical examples, we used the WENO scheme by using these kinds of mapping functions to calculate several discontinuous problems, such as one-dimensional SOD's problem, the interaction of shock-entropy wave, two-dimensional Riemann problem, implosion, and double Mach problem. Compared with the classical WENO-AIM, WENO- RM, and WENO-PM, these new adaptive mapped WENO schemes not only have higher resolution but also have high accuracy and computational stability.
\section*{Acknowledgement}
This work was supported by Scientific Research Foundation of Hunan Provincial Education Department (19C1766).


\begin{thebibliography}{99}
\setlength{\itemsep}{1pt}
\small{
\bibitem[1]{1} X. D. Liu, S. Osher, T. Chan, Weighted essentially non-oscillatory schemes, J. Comput. Phys. 115 (1994) 200-212.
\bibitem[2]{2} A. Harten, High resolution schemes for hyperbolic conservation laws, J. Comput. Phys. 49 (1983) 357-393.
\bibitem[3]{3} G.S. Jiang, C.W. Shu, Efficient implementation of weighted ENO schemes, J. Comput. Phys. 126 (1996) 202-228.
\bibitem[4]{4}D.S. Balsaraa, C. W. Shu. Monotonicity Preserving Weighted Essentially Non-oscillatory Schemes with Increasingly High Order of Accuracy. J. Comput.
Phys. 160 (2000) 405-452.
\bibitem[5]{5}J. Qiu, C. W. Shu. Hermite WENO schemes and their application as limiters for Runge-Kutta discontinuous Galerkin method: one-dimensional case. J. Comput. Phys. 193 (2003) 115-135.
\bibitem[6]{6}J. Qiu, C. W. Shu. Hermite WENO schemes and their application as limiters for Runge-Kutta discontinuous Galerkin method: two-dimensional case. Computers $\&$ Fluids. 34 (2005) 642-663.
\bibitem[7]{7}J. Qiu, C. W. Shu. Hermite WENO Schemes and Their Application as Limiters for Runge-Kutta Discontinuous Galerkin Method, III: Unstructured Meshes. J. Sci. Comput. 39 (2009) 293-321.
\bibitem[8]{8} A.K. Henrick, T.D. Aslam, J.M. Powers, Mapped weighted essentially non-oscillatory schemes: achieving optimal order near critical points, J. Comput.
Phys. 207 (2005) 542-567.
\bibitem[9]{9} R. Borges, M. Carmona, B. Costa, W.S. Don, An improved weighted essentially non-oscillatory scheme for hyperbolic conservation laws, J. Comput. Phys.
227 (2008) 3191-3211.
\bibitem[10]{10} H. Feng, F. Hu, R. Wang, A new mapped weighted essentially non-oscillatory scheme, J. Sci. Comput. 51 (2012) 449-473.
\bibitem[11]{11} H. Feng, C. Huang, R. Wang, An improved mapped weighted essentially non-oscillatory scheme, Appl. Math. Comput. 232 (2014) 453-468.
\bibitem[12]{12} R. Wang, H. Feng, C. Huang, A new mapped weighted essentially non-oscillatory method using rational mapping function, J. Sci. Comput. 67 (2016)
540-580.
\bibitem[13]{13} U. S. Vevek, B. Zang, New T.H., Adaptive mapping for high order WENO methods, J. Comput. Phys. 381 (2019) 162-188.
\bibitem[14]{14} Z. Hong, Z. Y. Ye, X. Z. Meng, A mapping-function-free WENO-M scheme with low computational cost, J. Comput. Phys. 405 (2020) 109145.
\bibitem[15]{15} F. X. Hu, High-order mapped WENO methods with improved efficiency, Computers $\&$ Fluids. 219 (2021) 104874.
\bibitem[16]{16} X. Wu, Y. Zhao. A high-resolution hybrid scheme for hyperbolic conservation laws, Int. J. Numer. Methods Fluids. 78 (2015) 162-187.
\bibitem[17]{17} P. Roe . Approximate Riemann Solvers, Parameter Vectors, and Difference Schemes, J. Comput. Phys. 43 (1981) 357-372.
\bibitem[18]{18} Y. Sun , Z. Wang. Evaluation of discontinuous galerkin and spectral volume methods for scalar and system conservation laws on unstructured grids. Int. J. Numer. Methods Fluids. 45(2004) 819-838.
\bibitem[19]{19} F. Davoudzadeh, H. McDonald, B. Thopson. Accuracy evaluation of unsteady CFD numerical schemes by vortex preservation. Computers $\&$ Fluids. 24 (1995) 883-895.
\bibitem[20]{20} S. Pirozzoli, On the spectral properties of shock-capturing schemes, J. Comput. Phys. 219 (2006) 489-497.
\bibitem[21]{21} R. Liska, B. Wendroff. Comparison of several difffference schemes on 1D and 2D test problems for the Euler equations. SIAM J.  Sci.  Comput. 25 (2003) 995-1017.}
\end{thebibliography}
\end{document}